\documentclass[12pt,a4paper]{amsart}

\usepackage{amsmath}\usepackage{amssymb}\usepackage{bbm}
\usepackage[T1]{fontenc}

\newcommand {\one} {{\bf 1}}

\usepackage{graphicx}
\usepackage{tabularx,mdwtab}
\usepackage{hyperref}
\usepackage{stmaryrd}

\usepackage{DLdef1}

\makeatletter
\theoremstyle{thm}
\@ifundefined{Theorem}{}{}
\@ifundefined{Lemma}{}{}
\@ifundefined{Corollary}{}{}
\@ifundefined{Statement}{}{}
\@ifundefined{Problem}{}{}
\@ifundefined{Question}{}{}
\theoremstyle{definition}
\@ifundefined{Remark}{}{}
\@ifundefined{Example}{}{}
\makeatother

\renewcommand {\ssbegin}[2][*]
 {\refstepcounter{subsection}%
\if#1*
\addcontentsline{toc}{subsection}{\thesubsection.\hskip 1pc #2}%
\else
\addcontentsline{toc}{subsection}{\thesubsection.\hskip 1pc #2. #1}%
\fi
 \def \secno {\gdef \secno {}{\ssecfont
\thesubsection.\hskip 2ex}%
 }%
 \begin{#2}}

\renewcommand {\sssbegin}[2][*]
 {\refstepcounter{subsubsection}%
\if#1*
\addcontentsline{toc}{subsubsection}{\thesubsubsection.\hskip 1pc #2}%
\else
\addcontentsline{toc}{subsubsection}{\thesubsubsection.\hskip 1pc #2. #1}
\fi
 \def \secno {\gdef \secno {}{\ssecfont \thesubsubsection.\hskip 2ex}%
 }%
 \begin{#2}}

\begin{document}

\title[Maximal subalgebras in three exceptional Lie superalgebras]{Maximal subalgebras
in three exceptional simple\\ finite-dimensional Lie superalgebras}

\author[LeLo]{Dimitry Leites${}^{a}$ and Oleksandr Lozhechnyk${}^{b}$}

\address{${}^{a}$Department of Mathematics, Stockholm University, Albanov\"agen 28, SE-114 19
Stockholm, Sweden; dimleites@gmail.com}

\address{${}^{b}$Faculty of Mechanics and Mathematics, Taras Shevchenko National University of
Kyiv, Hlushkova Avenue 4e, 03127 Kyiv, Ukraine; alozhechnik@gmail.com}

\date{2 October 2026}


\makeatletter
\@namedef{subjclassname@2020}{\textup{2020} Mathematics Subject Classification}
\makeatother

\begin{abstract}
Dynkin classified the maximal subalgebras of the simple finite-dimensional Lie algebras.
The super version of his problem is solved, for the matrix Lie superalgebras, by
Shchepochkina, and, for the simple vectorial ones with polynomial
coefficients taken with their Weisfeiler gradings, in the companion paper
by Leites and Shchepochkina --- with three exceptions left open: the three simple
finite-dimensional Lie superalgebras $\mathfrak{osp}_a(4|2)$, $\mathfrak{ag}(2)$ and $\mathfrak{ab}(3)$, which are
vectorial in more than one way and are not covered by the general argument.

Here we settle these three cases for \textbf{graded} subalgebras. For each of the three
Lie superalgebras, we list all the $\mathbb{Z}$-gradings to which the method of that paper
applies --- two for $\mathfrak{osp}_a(4|2)$, five for $\mathfrak{ag}(2)$, eight for $\mathfrak{ab}(3)$ --- and
describe, for every one of them, all the maximal graded subalgebras, the semi-simple ones
explicitly. A simple criterion upgrades most of the answers from ``maximal graded\rq\rq\ to
``maximal\rq\rq. Every answer is computed over $\mathbb{Q}$.

Two by-products: the classical list of \textbf{five} cases into which the graded
subalgebras of a vectorial Lie superalgebra split is incomplete for ambients of depth
$>1$ --- a sixth case has to be added, and it is non-empty precisely here; and the
parameter $a$ of the copies of $\mathfrak{osp}_a(4|2)$ sitting inside $\mathfrak{ag}(2)$,
resp. $\mathfrak{ab}(3)$, is computed intrinsically.
\end{abstract}

\maketitle

\markboth{\itshape Dimitry Leites\textup{,}  Oleksandr Lozhechnyk}{{\itshape Maximal subalgebras in three exceptional Lie superalgebras}}

\thispagestyle{empty}
\setcounter{tocdepth}{2}
\tableofcontents

\section{Introduction}\label{Sintro}

\ssec{Notation}\label{ssNotation} We work over $\Cee$; all the models are, however,
defined over $\Qee$, and so is every computation reported below.

$\sdim V:=(\dim V_\ev|\dim V_\od)$ is the superdimension of a~superspace $V$; the
superdimension is given for every superspace and every Lie superalgebra we name, and the
rank is given for every operator we name. Further:

$\id_\fg$ is the tautological module over a~linear Lie (super)algebra $\fg$.

$\pi_i$ is the $i$th fundamental weight of a~simple Lie algebra, the simple roots being
numbered as in \cite{OV}; the same symbol denotes the irreducible module with highest
weight $\pi_i$, so that $\pi_1$ over $\fg(2)$ is its $7$-dimensional module.

$\spin_n$ is the spinor module of $\fo(n)$, so $\dim\spin_7=\dim\spin_8=8$. 

$\Pi$ is the
change-of-parity functor.

$\fc(\fg):=\fg\oplus\Cee z$ is the trivial central extension of $\fg$ (as in
\cite{LSh}); for example, $\fc(\fosp(2|2))\simeq\fgl(1|2)$.

$\fst(V):=\{X\in\fg\mid XV\subseteq V\}$ is the stabilizer of a~subspace~$V$.

$\fhei(2n|m)$ is the Heisenberg Lie superalgebra, $\sdim\fhei(2n|m)=(2n+1|m)$.

$\fk(2n+1|m)$ is the contact Lie superalgebra.

$\ft$ is a~maximal toral subalgebra of the ambient $\fg$; the roots and the weights are
written in the coordinates fixed in tables \eqref{tabAgGr}, \eqref{tabAbGr}.

If $\fg_0=\fa_1\oplus\dots\oplus\fa_k$ is a~direct sum of algebras and $V_i$ is
an $\fa_i$-module for each $i$, then\index{$\boxtimes$}
\begin{equation}\label{outer}
V_1\boxtimes\dots\boxtimes V_k
\end{equation}
denotes the \textit{outer} tensor product: the superspace
$V_1\otimes\dots\otimes V_k$ on which each summand $\fa_i$ acts through the
$i$th factor alone, \textit{i.e.},
\[
\text{$a(v_1\otimes\dots\otimes
v_k)=(-1)^{p(a)(p(v_1)+\dots+p(v_{i-1}))}v_1\otimes\dots\otimes av_i
\otimes\dots\otimes v_k$ for every $a\in\fa_i$.}
\] 

Following \cite{LSh}, the adjective ``simple\rq\rq\ is, by abuse of notation, broadened so as
to include also the ``relatives\rq\rq\ of simple algebras (central extensions, algebras of
derivations, and so on); ``almost simple\rq\rq\ means that the algebra sits between a~simple
one and its algebra of derivations.

\ssec{Description of the ambients} Following \cite{Ld}, we do not use the early notation ``\`a~la~Cartan\rq\rq\ for the
simple Lie superalgebras: we write $\fosp_a(4|2)$ for the members of the one-parameter
family sometimes denoted $D(2,1;a)$, use  Serganova's suggestive
notation $\fag(2)$ and $\fab(3)$, see \cite{Se4}, for the two exceptions sometimes denoted $G(3)$ and $F(4)$, respectively, whose even part is the sum of $\fsl(2)$ and
$\fg(2)$, resp.~$\fb_3=\fo(7)$. Cartan's notation for the root
systems of Lie \textbf{algebras} ($A_1$, $B_3$, $G_2$) is retained.

\sssec{Explicit models of $\fag(2)$ and $\fab(3)$}\label{Smodels}
Both these superalgebras are of the shape
\be\label{shape3exc}
\fg=(\fs\oplus\fsl(2))\oplus(V\boxtimes\Cee^2),\qquad
\begin{cases}
\fs=\fg(2),\ V=\pi_1\ (\dim=7)&\text{for $\fag(2)$},\\
\fs=\fo(7),\ V=\spin_7\ (\dim=8)&\text{for $\fab(3)$},
\end{cases}
\ee
the even part acting in the obvious way, the odd bracket being
\be\label{oddbr3exc}
[v\otimes u,v'\otimes u']=\alpha\,(v,v')_V\,\gamma(u,u')+\beta\,\omega(u,u')\,
\pi(v\wedge v'),
\ee
where $(- ,- )_V$ is the non-degenerate invariant symmetric form on $V$, $\omega$ is the
symplectic form on $\Cee^2$, and ${\gamma\colon S^2(\Cee^2)\tto\fsl(2)}$ is an
isomorphism of $\fsl(2)$-modules (unique up to a~factor, which is fixed in our models
once and for all), and
\be\label{moment3exc}
\pi\colon\Lambda^2V\tto\fs,\qquad
\pi(v\wedge v')=\sum_i\langle A_iv,v'\rangle_V\,A_i^*
\ee
is the moment map, where $A_i^*$ is the basis vector of $\fs$ dual to  $A_i$ with respect
to the trace form. 

(Note that the shape \eqref{oddbr3exc} is forced: since
$S^2(V\otimes\Cee^2)=S^2V\otimes S^2\Cee^2\oplus\Lambda^2V\otimes\Lambda^2\Cee^2$,
and since (in both cases) $\fs\not\subset S^2V$, whereas
$\Lambda^2V=\fs\oplus V'$ and $S^2V\supset\Cee$, the space of
$\fg_\ev$-morphisms ${S^2(\fg_\od)\tto\fg_\ev}$ is exactly $2$-dimensional. The two
constants are determined up to the rescaling $v\otimes u\mapsto t(v\otimes u)$ of
$\fg_\od$ by a~number $t\neq0$, which multiplies $\alpha$ and $\beta$ by $t^2$,
by the %$(\od,\od,\od)$ 
Jacobi identity for 3 odd elements; the identity is linear in $\alpha$ and in $\beta$
separately; solving it gives a~$1$-dimensional space of solutions:
\[
\fag(2):\ \ \tfrac{\beta}{\alpha}=3,\qquad\qquad \fab(3):\ \ \tfrac{\beta}{\alpha}=\tfrac{8}{3};
\]
the ratio refers to the normalizations of $\gamma$ and of the trace form fixed in our
models, and changes together with them, whereas its non-vanishing is intrinsic.)
\iffalse for the normalizations fixed below. Both models are realized over $\Qee$:

1) $\fg(2)$ is realized as the stabilizer, inside $\fgl(7)$, of the $3$-form
$\varphi$ of weight $0$ in the weight basis $w_0;w_1,w_2,w_3;w_{-1},w_{-2},w_{-3}$
of the $7$-dimensional module (of weights $0$, $\pm\eps_i$ with
$\eps_1+\eps_2+\eps_3=0$),
\[
\varphi=w_1\wedge w_2\wedge w_3+w_{-1}\wedge w_{-2}\wedge w_{-3}+
\textstyle\sum_iw_0\wedge w_i\wedge w_{-i}
\]
(any choice of non-zero constants in this ansatz gives a~$14$-dimensional
stabilizer.%\AO{DL: verified by computation.}). 
Its weights are the $12$ roots
of $G_2$ plus $0$ twice, as they should be, and the invariant symmetric form on
$V$ is recovered from $\fg(2)$, not postulated.

2) $\spin_7$ is realized on $\Lambda^\bullet(W)$, where $W=\Span(f_1,f_2,f_3)$, by
the Clifford action normalized by $c(x)c(y)+c(y)c(x)=(x,y)\Id$ (with
$(w_0,w_0)=2$, which is what makes $c(w_0)=(-1)^{\deg}$ rational), and
$\fo(7)\subset\fgl(\spin_7)$ is recovered from the relation $[\rho(A),c(z)]=c(Az)$;
that $\rho$ is a~homomorphism, and the invariant symmetric form on $\spin_7$, are
computed, not postulated.
\fi

\ssec{Description of $\fg=\fosp_a(4|2)$}\label{ssOspModel}
Explicitly, for $\fg_\ev=\fsl(2)^{\oplus3}$ and
$\fg_\od=\Cee^2\boxtimes\Cee^2\boxtimes\Cee^2$ with
\be\label{modelOsp}
[x_1\otimes x_2\otimes x_3,\ y_1\otimes y_2\otimes y_3]=
\sum_{1\leq m\leq 3}c_m\Big(\prod_{k\neq m}\omega(x_k,y_k)\Big)\gamma_m(x_m,y_m),
\ee
where $\omega$ is the symplectic form on $\Cee^2$ and
$\gamma\colon S^2(\Cee^2)\tto\fsl(2)$ is the standard isomorphism; the parameter
is $(c_1:c_2:c_3)$ and $a:=\tfrac{c_2}{c_1}$ for $c=(1,a,-1-a)$; for the normalization of $a$ and the
six-element group acting on it, see \cite{Ld}. The answers below are given in the
normalization-free form \eqref{jinv} as well, so that nothing depends on that choice. (Kaplansky, who discovered the three exceptional Lie superalgebras, proved that the Jacobi identity is satisfied if and only if $c_1+c_2+c_3=0$, see \cite{Ld}.)

\textbf{Normalization of the Cartan matrix}. The row of a~Cartan matrix $A$ that
corresponds to an~isotropic simple root $\alpha_i$ \textup{(}the one with $A_{ii}=0$\textup{)}
is only defined up to a~non-zero factor. We fix this factor by requiring the
entry of that row in the \textbf{left} column to be equal to~$-1$. For $\fosp_a(4|2)$
realized with two non-isotropic simple roots, this yields
\be\label{cmOsp}
\begin{pmatrix} 2&-1&0\\ -1&0&-a\\0&-1&2\end{pmatrix},
\ee
where the parameter is precisely the $a=\tfrac{c_2}{c_1}$ of eq.~\eqref{modelOsp}. %; this is the matrix used in \cite{LS1}. 
The all-odd system of simple roots we keep, as usual, in the
symmetric \textup{(}Gram\textup{)} form, since the rule applied to it literally introduces
denominators: it turns, e.g., $(1+a, -a, 0)$ into $(-1, \tfrac{a}{1+a},0)$, destroying the
shape in which such matrices are always printed. % --- the same happens to the matrix of $\fsvect^L_a(1|2)$ in \cite[p.~188]{Ber}. 
So the normalization
applies to the matrices with two non-isotropic simple roots, and the Cartan matrices
of the all-odd systems of simple roots stay symmetric. 

\sssec{The parameter $a$ and $j(a)$} The isomorphisms $\fosp_{\alpha}(4|2)\simeq \fosp_{\alpha'}(4|2)$ are generated by the
transformations (see \cite{BGL}):
\begin{equation}\label{osp42symm}
\alpha\longmapsto \alpha':= -1-\alpha\ , \qquad \alpha\longmapsto \alpha':=
\nfrac{1}{\alpha}\ ,
\end{equation}
so the other isomorphisms are 
\begin{equation}\label{osp42symm1}
\alpha\longmapsto \alpha':=
-\nfrac{1+\alpha}{\alpha}\ ,\quad\alpha\longmapsto \alpha':=
-\nfrac{1}{\alpha+1}\ ,\quad\alpha\longmapsto \alpha':= -\nfrac{\alpha}{\alpha+1}.
\end{equation}

The parameter $a$ is defined up to the action of the group $S_3$ generated by the
transformations \eqref{osp42symm}, so a~value of $a$ is best given in terms of the
following $S_3$-invariant,
\be\label{jinv}
j(a):=\frac{(a^2+a+1)^3}{a^2(a+1)^2},
\ee
which does not depend on any normalization of the parameter.

\ssec{The maximal subalgebras in the three exceptional simple ambients}\label{ssProblem} S.~Lie was the first to try to classify the maximal
subalgebras of the simple vectorial (i.e., consisting of vector fields) Lie algebras with
polynomial coefficients, see \cite{Lie}. As stated, the problem is wild, as proved in  \cite{ShM}. Dynkin
\cite{Dy} solved the tame problem: he classified the maximal subalgebras
of the simple finite-dimensional Lie algebras. Shchepochkina solved the version of Dynkin's problem for
\textbf{matrix} Lie superalgebras, see \cite{ShM}, where,
however, three ambients were left aside:
\be\label{3exc}
\begin{array}{c}
\fosp_a(4|2),\qquad \fag(2),\qquad
\fab(3).
\end{array}
\ee 

For the simple
\textbf{vectorial} Lie superalgebras with polynomial coefficients, considered with the
gradings associated with their Weisfeiler filtrations (W-filtrations for short), the problem is solved in \cite{LSh}. 
The three Lie superalgebras \eqref{3exc} are simple,
finite-dimensional, and \textbf{vectorial in several ways at once}: each of them admits
several Weisfeiler $\Zee$-gradings, and hence several realizations by vector
fields,  see \cite{Sh14, LSh};
the arguments of \cite{LSh} are designed for infinite-dimensional ambients of any
depth, and for finite-dimensional ones only in depth~$1$. So, of the three exceptions,
they cover only $\fosp_a(4|2)$ and $\fab(3)$, and only in their gradings of depth~$1$,
see Subsec.~\ref{ssTKK}; they cover no grading of $\fag(2)$, and none of depth $2$
or~$3$. Filling this gap is the purpose of this paper. 

The three exceptions \eqref{3exc} are also realized as the symmetry superalgebras of
supergeometries on the flag supermanifolds $G/P$, for every parabolic subgroup $P$:
$\fag(2)$ and $\fab(3)$ in \cite{KST, SaT, KLl}, and $\fosp_a(4|2)$ in \cite{EKT}.
The $\Zee$-gradings behind these realizations are the ones defined by all the parabolic
subalgebras, and the Cartan prolong of the negative part is computed there for each of
them. Here only the gradings satisfying the conditions \eqref{cond3exc} are needed; how
they sit among the parabolic ones is described in Remark~\ref{RParab3exc}.

The even and the odd parts of the three exceptions \eqref{3exc} are
\be\label{3excev}
\begin{array}{lll}
\fg=\fosp_a(4|2):&\fg_\ev=\fsl(2)\oplus\fsl(2)\oplus\fsl(2),&
\fg_\od=\Cee^2\boxtimes\Cee^2\boxtimes\Cee^2\ (\sdim=0|8),\\
\fg=\fag(2):&\fg_\ev=\fg(2)\oplus\fsl(2),&
\fg_\od=\pi_1\boxtimes\Cee^2\ (\sdim=0|14),\\
\fg=\fab(3):&\fg_\ev=\fo(7)\oplus\fsl(2),&
\fg_\od=\spin_7\boxtimes\Cee^2\ (\sdim=0|16).
\end{array}
\ee

\ssec{What is done here}\label{ssPlan} For each of the three exceptions, we

\textup{(i)} list the $\Zee$-gradings to which the method of \cite{LSh} applies, see
\S\ref{ss3excGrad}; there are two of them for $\fosp_a(4|2)$, five for $\fag(2)$ and
eight for $\fab(3)$, all of depth $\leq3$;

\textup{(ii)} describe, for every such grading, all the maximal \textbf{graded}
subalgebras, see \S\ref{Sosp} for $\fosp_a(4|2)$ and \S\ref{Sag}--\S\ref{Sab} for the two
other exceptions.

The answers marked \textbf{max} are maximal subalgebras of $\fg$, not only maximal among
the graded ones, thanks to the criterion of \S\ref{ss3excCert}; the answers marked
\textbf{max(b)} are maximal graded, and no intermediate subalgebra was found, but the
verification is not exhaustive, see Remark \ref{Rb3exc}.

For $\fosp_a(4|2)$ everything is done by hand, see Theorem \ref{ThOsp3exc}. For $\fag(2)$
and $\fab(3)$ we use explicit models over $\Qee$, see \S\ref{Smodels}.

\ssec{The answers}\label{ssAnswers} Collected over all the gradings, the semi-simple
maximal graded subalgebras of the three exceptions are the following. Recall that
$\fg_\ev$ is one of them for each of the three, and for every grading, see Lemma
\ref{L3-3exc}.

\sssbegin[The semi-simple maximal graded subalgebras]{Theorem}[The semi-simple maximal graded subalgebras]\label{ThMain} \textup{1)} For $\fg=\fosp_a(4|2)$, the semi-simple maximal
graded subalgebras are, up to conjugacy, exactly
\[
\begin{array}{c}
\fsl(2)\oplus\fsl(2)\oplus\fsl(2)=\fg_\ev,\quad
\fsl(2)\oplus\fosp(1|2),\ \ \ 
\fosp(3|2),
\end{array}
\]
the first two for every $a$, the third one only for $a\in\{1,-2,-\nfrac12\}$, i.e.,
precisely when $\fosp_a(4|2)\simeq\fosp(4|2)$.

\textup{2)} For $\fg=\fag(2)$, they are exactly the seven of the list \eqref{unionAg}, namely
\[
\begin{array}{lll}
\fg(2)\oplus\fsl(2)=\fg_\ev\ (17|0),&\fosp_a(4|2)\ \text{with $j(a)=\nfrac{13^3}{144}$}\
(9|8),&\fosp(3|2)\oplus\fsl(2)\ (9|6),\\
\fsl(3|1)\ (9|6),&\fsl(3)\oplus\fosp(1|2)\ (11|2),&\fsl(2)\oplus\fsl(2|1)\ (7|4),\\
\fosp(1|2)\oplus\fosp(1|2)\ (6|4).&&
\end{array}
\]

\textup{3)} For $\fg=\fab(3)$, they are exactly the six of the list \eqref{unionAb}, namely
\[
\begin{array}{lll}
\fo(7)\oplus\fsl(2)=\fg_\ev\ (24|0),&\fsl(4|1)\ (16|8),&
\fosp_a(4|2)\oplus\fsl(2)\ \text{with $j(a)=\nfrac{7^3}{36}$}\ (12|8),\\
\fg(2)\oplus\fosp(1|2)\ (17|2),&\fsl(2|1)\oplus\fsl(3)\ (12|4),&
\fosp(1|2)\oplus\fosp(1|2)\oplus\fsl(2)\ (9|4).
\end{array}
\]

In headings \textup{2)} and \textup{3)}, $j$ is the invariant \eqref{jinv} of the
six-element group acting on the parameter $a$, and $\sdim$ is given in parentheses. Every
maximal graded subalgebra other than these is either $\fg_{\leq0}$, or $\fg_{\geq0}$, or
the stabilizer-with-prolongation of a~subspace of $\fg_{-1}$, i.e., a~parabolic-like
answer, cf.~Remark $\ref{Rthin3exc}$.
\end{Theorem}

The proof of Theorem \ref{ThMain} is the content of \S\S\ref{Sosp}--\ref{Sab}: heading 1)
is Theorem \ref{ThOsp3exc}, heading 2) is Statements \ref{SAg1}--\ref{SAg5}, and heading
3) is Statements \ref{SAb1}--\ref{SAb8}, where, moreover, the answers are given grading by
grading, with their graded components.

\ssec{Two by-products}\label{ssByProducts} \textup{1)} The set of graded subalgebras of
a~$\Zee$-graded vectorial Lie superalgebra splits, in \cite{LSh} and in its predecessors,
into five non-intersecting cases, see the list of cases \eqref{6cases}. For ambients of depth $d>1$, these
five cases are \textbf{not} exhaustive: the case $\fu_{-1}=0\neq\fu_-$ is missing. It is
non-empty, and the examples that show this are exactly the three exceptions in their
contact gradings, see Statement \ref{6o3exc}.

\textup{2)} The copies of $\fosp_a(4|2)$ that occur in the answers have a~perfectly
definite value of the parameter $a$, and this value is computed here intrinsically, in
a~form independent of any normalization of $a$, see Remark \ref{Ra3exc}. The reason for
the two values that occur is the ratio of the lengths of the roots in $G_2$, resp.~$B_3$.

\section{Preliminaries}\label{Sprelim}

Everything in this section is taken from \cite{LSh}; we reproduce it in order to make the
paper self-contained.

\ssec{$\Zee$-gradings, Weisfeiler gradings, depth}\label{ssGradings} Let
$\fg=\oplus_{i\geq-d}\ \fg_i$ be a~$\Zee$-graded Lie superalgebra with $\fg_{-d}\neq0$;
the number $d$ is the \textit{depth} of the grading. The grading is said to be
\textit{without gaps} if $\fg_i\neq0$ for all $-d\leq i\leq d'$, where $d'$ is the highest
degree, and \textit{symmetric} if $\dim\fg_i=\dim\fg_{-i}$ for all $i$ --- which is
automatic if $\fg$ is finite-dimensional with a~non-degenerate invariant symmetric bilinear
form, as is the case for the three exceptions \eqref{3exc}.

The grading is a~\textit{Weisfeiler grading} (\textit{W-grading} for short) if
$\fg_{-1}$ generates $\fg_-:=\oplus_{i<0}\ \fg_i$, if the $\fg_0$-action on $\fg_{-1}$ is
irreducible, and if the grading is transitive, i.e., if
$\{x\in\fg_k\mid[x,\fg_-]=0\}=0$ for $k\geq0$; the W-gradings of the simple vectorial Lie
superalgebras with polynomial coefficients are classified in \cite{LS1}, and the depth of
such a~grading is 1, 2 or~3.

\ssec{The Cartan prolongs}\label{ssProlong} Let $\fg_-=\oplus_{-d\leq i\leq-1}\ \fg_i$ be
a~nilpotent $\Zee$-graded Lie superalgebra, and let $\fg_0$ be a~subalgebra of the Lie
superalgebra $\fder_0(\fg_-)$ of the degree-$0$ derivations of~$\fg_-$. The maximal
transitive $\Zee$-graded Lie superalgebra $\fg=\oplus_{k\geq-d}\ \fg_k$ whose non-positive
part is the given $\fg_-\oplus\fg_0$ is called the \textit{Cartan prolong} of the pair
$(\fg_-,\fg_0)$ and is denoted by $(\fg_-,\fg_0)_*$. Explicitly, and recursively, one degree at a~time,
\be\label{prol3exc}
\fg_k=\{x\mid[x,\fg_j]\subseteq\fg_{j+k}\text{~~for all $j<0$}\}\qquad\text{for $k\geq1$}.
\ee
Given a~$\Zee$-graded ambient $\fg$ and a~graded subalgebra $\fu$ of it with a~prescribed
non-positive part $\fu_{\leq0}$, the largest graded subalgebra of $\fg$ with that
non-positive part is obtained by the same recursion \eqref{prol3exc} performed inside
$\fg$; this is the \textit{Cartan prolong of $\fu_{\leq0}$ in $\fg$}, and it is what makes
the search for maximal graded subalgebras finite.

\ssec{The problem, restricted to graded subalgebras}\label{ssPROB} We consider only \textbf{graded} subalgebras; this is
a~vital restriction. The problem, as stated in \cite{LSh}, is:
\be\label{PROB}
\begin{minipage}[l]{14cm}
\textbf{In ``simple\rq\rq\ vectorial Lie superalgebras $\fg=\oplus_{i\geq-d}\ \fg_i$, such as
Cartan prolongs of $\fg_-$, classify maximal simple, semi-simple, almost simple and
semi-almost simple graded Lie subsuperalgebras
$\fu=\oplus_{i\geq-d}\ (\fu\cap\fg_i)$}.
\end{minipage}
\ee
Here, and everywhere below, $\fu_i:=\fu\cap\fg_i$.

\ssec{The six cases}\label{ssCases} Whatever the depth of $\fg$, the set of graded
subalgebras $\fu\subset\fg$, where $\fg$ is either simple or a~Cartan prolong of $\fg_-$
with an irreducible $\fg_0$-action on $\fg_{-1}$, splits into the following
non-intersecting subsets:
\begin{equation}\label{6cases}
\begin{minipage}[l]{14cm}
$1^o$. $\fu_{-1}=\fg_{-1}$ and $\fu_{0}=\fg_{0}$; $\fu_{1}\neq 0$;\\
$2^o$. $0 \subsetneq \fu_{-1}\subsetneq \fg_{-1}$;\\
$3^o$. $\fu_{-1}=\fg_{-1}$ and $\fu_{0}\neq \fg_{0}$;\\
$4^o$. $\fu_{-1}=\fg_{-1}$ and $\fu_{0}=\fg_{0}$; $\fu_{1}=0$;\\
$5^o$. $\fu_{-}=0$, where $\fu_-= \oplus_{i<0}\ (\fu\cap\fg_{i})$;\\
$6^o$. $\fu_{-1}=0$, but $\fu_-\neq0$.
\end{minipage}
\end{equation}
For ambients of depth $d=1$ the case $6^o$ is empty, and \eqref{6cases} is the list of
five cases used in \cite{LSh}. For $d>1$, the case $6^o$ is not empty; the examples are
given in Statement \ref{6o3exc}, and they are what forced the sixth item.

\ssbegin[Obvious subcases of the cases \eqref{6cases}]{Lemma}[Obvious subcases of
the cases \eqref{6cases}]\label{Case 4} In case \textup{$4^o$}, $\fu=\oplus_{i\leq0}\ \fg_i$ is
maximal in~$\fg$. In case \textup{$5^o$}, $\fu$ is a~maximal subalgebra of $\fg$ if and
only if $\fu=\oplus_{i\geq0}\ \fg_i$. \end{Lemma}

\begin{proof} The irreducibility of the $\fg_0$-action on $\fg_{-1}$ implies both
statements. \end{proof}

\ssec{A criterion of maximality}\label{ss3excCert} The following simple observation is
what turns the lists below into certificates rather than samples; it applies to any
ambient with $1$-dimensional root spaces, in particular to each of the three superalgebras \eqref{3exc}, and it
strengthens the conclusion from ``maximal graded\rq\rq\ to ``maximal\rq\rq.

\ssbegin[What upgrades maximal graded to maximal]{Lemma}[What upgrades maximal graded to maximal]\label{Lcert3exc}
Let $\fg$ be a~Lie superalgebra with a~maximal toral subalgebra $\ft$ such that all the
root spaces are $1$-dimensional and all the roots are distinct \textup{(}this is
the case for each of the three superalgebras \eqref{3exc}\textup{)}. Let $\fu\subset\fg$ be
a~subalgebra and $\ft'=\fu\cap\ft$. Then, every subalgebra $\fw$ such that
$\fu\subset\fw\subset\fg$ is $\ft'$-invariant, and hence is the sum of its
intersections with the $\ft'$-weight spaces of~$\fg$. In particular:

\textup{1)} if $\ft\subset\fu$, then all the weight spaces of $\fg/\fu$ are
$1$-dimensional, and $\fu$ is a~maximal subalgebra~of $\fg$ if and only if the
subalgebra generated by $\fu$ and $x_\alpha$ is equal to $\fg$ for every root
vector $x_\alpha\not\in\fu$ --- a~finite verification;

\textup{2)} if $\ft\not\subset\fu$, the same test applied to the $1$-dimensional
$\ft'$-weight spaces of $\fg/\fu$ is still conclusive, and only the weight spaces
of dimension $>1$ have to be sampled.
\end{Lemma}

\begin{proof} Since $\ft'\subset\fu\subset\fw$ and $\fw$ is a~subalgebra, we have
$[\ft',\fw]\subset\fw$; as $\ft'$ acts on $\fg$ semi-simply, its invariant
subspace $\fw$ is the sum of its weight components. If $\ft\subset\fu$, then the
weight-$0$ component of $\fg$ is $\ft$ itself, and hence lies in $\fu$, whereas
the other weight spaces are the root spaces, of dimension $1$; therefore
a~$\fw\supsetneq\fu$ contains a~whole root line $\Cee x_\alpha$ with
$x_\alpha\not\in\fu$, and hence
$\fw\supset\langle \fu,x_\alpha\rangle=\fg$. \end{proof}

In the tables below, the answers obtained with $\ft\subset\fu$ are marked \textbf{max}:
they are maximal subalgebras of the ambient. The ones for which some $\ft'$-weight space
of $\fg/\fu$ is of dimension $>1$ are marked \textbf{max(b)}: there the test was run on
the basis vectors of these weight spaces and on random vectors in them, so such an entry
means ``no intermediate subalgebra found\rq\rq, and is not a~proof; see Remark \ref{Rb3exc}.

\section{Step \textup{(i)}: the gradings to which the method applies}\label{ss3excGrad}

The method of \cite{LSh} applies to a~$\Zee$-graded $\fg=\oplus_{i\geq-d}\ \fg_i$ such
that
\be\label{cond3exc}
\begin{minipage}[l]{14cm}
$\bullet$ the grading is symmetric, of depth $d$, without gaps;\\
$\bullet$ $\fg_{-1}$ generates $\fg_-$ \textup{(}the Weisfeiler
property\textup{)};\\
$\bullet$ the $\fg_0$-action on each $\fg_{-k}$, where $k>0$, is irreducible.
\end{minipage}
\ee
Every $\Zee$-grading of $\fg$ is given by an element $h$ of a~maximal toral subalgebra
such that $\alpha(h)\in\Zee$ for every root $\alpha$, the degree of $\fg_\alpha$
being~$\alpha(h)$. Since all the root spaces are $1$-dimensional, and $\fg_{-k}$
consists of root vectors only for $k>0$, every $\fg_0$-submodule of $\fg_{-k}$ is
spanned by root vectors, and therefore
\be\label{irr3exc}
\text{$\fg_{-k}$ is $\fg_0$-irreducible}\iff\text{the set of roots of $\fg_{-k}$
is connected under adding roots of degree $0$,}
\ee
which turns the classification of the gradings satisfying \eqref{cond3exc} into
a~finite check on the root system. Its result is as follows.

\ssbegin[The gradings to which the method applies: $2+5+8$]{Statement}[The gradings to which the method applies: $2+5+8$]\label{StGr3exc} Up to isomorphism, the gradings of the three
exceptions \eqref{3exc} satisfying \eqref{cond3exc} are exactly the ones listed
in tables \eqref{tabOspGr}--\eqref{tabAbGr}: two for $\fosp_a(4|2)$, five for
$\fag(2)$, eight for $\fab(3)$; all of them are of depth~$\leq3$. If depth $4$ is
allowed, exactly one more grading appears for $\fag(2)$ and one for $\fab(3)$,
namely
\[
\begin{array}{ll}
\fag(2):&1|1,\ 2|2,\ 1|1,\ 2|2,\ \mathbf{5|2},\ 2|2,\ 1|1,\ 2|2,\ 1|1,\\
\fab(3):&3|0,\ 0|2,\ 3|0,\ 0|6,\ \mathbf{12|0},\ 0|6,\ 3|0,\ 0|2,\ 3|0
\end{array}
\]
\textup{(}the superdimensions of the graded components, the $0$th one in
boldface\textup{)}; there are no gradings of depth~$\geq5$. The grading elements $h$
were enumerated over the weight lattice, and \eqref{irr3exc} applied to each.
\end{Statement}

\be\label{tabOspGr}
\begin{array}{c}
\text{$\fg=\fosp_a(4|2)$, $\sdim \fg=9|8$; \ $h$ is given by $(\eps_1(h),\eps_2(h),\eps_3(h))$}\\
\renewcommand{\arraystretch}{1.35}
\begin{tabular}{|c|c|c|c|c|c|}
\hline
$N$&$h$&$d$&$\sdim\fg_{-d},\dots,\sdim\fg_0$&$\fg_0$&ambient\cr
\hline\hline
1&$(\nfrac12,\nfrac12,0)$&1&$2|2;\ 5|4$&$\fc(\fosp(2|2))$&
$(\fg_{-1},\fg_0)_*$\cr
\hline
2&$(1,0,0)$&2&$1|0,\ 0|4;\ 7|0$&$\fc(\fo(4))=\fk(1|4)_0$&$\fk(1|4)$\cr
\hline
\end{tabular}
\end{array}
\ee

\be\label{tabAgGr}
\begin{array}{c}
\text{$\fg=\fag(2)$, $\sdim \fg=17|14$; \ $h$ is given wrt
$(\delta;\eps_1,\eps_2,\eps_3)$}\\
\renewcommand{\arraystretch}{1.35}
\begin{tabular}{|c|c|c|c|c|c|}
\hline
$N$&$h$&$d$&$\sdim\fg_{-d},\dots,\sdim\fg_0$&$\fg_0$&ambient\cr
\hline\hline
1&$(1;0,0,0)$&2&$1|0,\ 0|7;\ 15|0$&$\fc(\fg(2))$&$\fk(1|7)$\cr
\hline
2&$(0;0,1,-1)$&2&$1|0,\ 4|4;\ 7|6$&$\fc(\fosp(3|2))$&$\fk(5|4)$\cr
\hline
3&$(1;0,1,-1)$&2&$2|2,\ 4|3;\ 5|4$&$\ast$&$(\fg_-,\fg_0)_*$\cr
\hline
4&$(0;1,1,-2)$&3&$2|0,\ 1|2,\ 2|4;\ 7|2$&$\ast$&$(\fg_-,\fg_0)_*$\cr
\hline
5&$(1;1,1,-2)$&3&$2|1,\ 2|2,\ 2|2;\ 5|4$&$\ast$&$(\fg_-,\fg_0)_*$\cr
\hline
\end{tabular}
\end{array}
\ee

\be\label{tabAbGr}
\begin{array}{c}
\text{$\fg=\fab(3)$, $\sdim \fg=24|16$; \ $h$ is given wrt
$(\delta;\eps_1,\eps_2,\eps_3)$}\\
\renewcommand{\arraystretch}{1.35}
\begin{tabular}{|c|c|c|c|c|c|}
\hline
$N$&$h$&$d$&$\sdim\fg_{-d},\dots,\sdim\fg_0$&$\fg_0$&ambient\cr
\hline\hline
1&$(1;0,0,1)$&1&$6|4;\ 12|8$&$\fc(\fosp(2|4))$&$(\fg_{-1},\fg_0)_*$\cr
\hline
2&$(2;0,0,0)$&2&$1|0,\ 0|8;\ 22|0$&$\fc(\fo(7))$&$\fk(1|8)$\cr
\hline
3&$(0;0,1,1)$&2&$1|0,\ 6|4;\ 10|8$&$\ast$&$\fk(7|4)$\cr
\hline
4&$(0;0,0,2)$&2&$5|0,\ 0|8;\ 14|0$&$\fsl(2)\oplus\fsp(4)\oplus\Cee$&
$(\fg_-,\fg_0)_*$\cr
\hline
5&$(1;1,1,1)$&2&$3|1,\ 4|4;\ 10|6$&$\ast$&$(\fg_-,\fg_0)_*$\cr
\hline
6&$(2;0,1,1)$&2&$2|2,\ 6|4;\ 8|4$&$\ast$&$(\fg_-,\fg_0)_*$\cr
\hline
7&$(0;1,1,2)$&3&$2|0,\ 2|2,\ 4|4;\ 8|4$&$\ast$&$(\fg_-,\fg_0)_*$\cr
\hline
8&$(3;1,1,1)$&3&$1|1,\ 3|3,\ 3|3;\ 10|2$&$\ast$&$(\fg_-,\fg_0)_*$\cr
\hline
\end{tabular}
\end{array}
\ee

In table \eqref{tabOspGr}, $\pm2\eps_i$ are the roots of the $i$th summand $\fsl(2)$ of
$\fg_\ev$ and $\pm\eps_1\pm\eps_2\pm\eps_3$ are the odd roots, so the degree of
$\fg_\alpha$ is $\alpha(h)$. Here, $\ast$ marks the components not identified in intrinsic terms (their even
parts and the $\fg_0$-module structure of $\fg_{-1}$ are, of course, known from
the root data); ``ambient\rq\rq\ is the vectorial Lie superalgebra in which the grading
realizes $\fg$, i.e., the one containing the Cartan prolong of $(\fg_-,\fg_0)$.

\ssec{Depth $1$ occurs only twice, and it is the Tits--Kantor--Koecher
grading}\label{ssTKK} The Lie superalgebra $\fosp_a(4|2)$ has exactly one grading of
depth $1$, with $\fg_{-1}$ of $\sdim=(2|2)$, and $\fab(3)$ has exactly one, with $\fg_{-1}$
of $\sdim=(6|4)$, whereas $\fag(2)$ has \textbf{none}. This matches the
Tits--Kantor--Koecher picture: the $\Zee$-gradings of depth $1$ correspond to
Jordan superalgebras, see \cite{K77, Kap}, and $\fg_{-1}$ is, respectively, the
$4$-dimensional Jordan superalgebra $D_t$ of the $1$-parameter family (whence the
parameter $a$) and the $10$-dimensional Kac Jordan superalgebra $K_{10}$, whose even part
is the Jordan algebra of a~non-degenerate quadratic form on $\Cee^{1+5}$ and whose odd
part is the spinor module $\spin_5$ of $\fo(5)=\fsp(4)$, i.e., the tautological
$\fsp(4)$-module, $\dim\spin_5=4$, --- exactly what the root data of
table \eqref{tabAbGr} give. For $\fag(2)$, the corresponding Kaplansky superalgebra $K_3$,
$\sdim K_3=(1|2)$, is not unital and, accordingly, no grading of depth $1$ exists; the
shallowest grading of $\fag(2)$ is the contact one, $N=1$ in table \eqref{tabAgGr}.

Consequently, of the three exceptions the complete answer given in \cite{LSh} for
depth-$1$ ambients covers only $\fosp_a(4|2)$ and $\fab(3)$; for $\fag(2)$ one has
to work in depth $2$ from the start.

\ssec{The uniform contact family}\label{ssContactFamily} Each of the three exceptions has
a~contact grading, with $\fg_{-2}=\Cee$ even and $\fg_{-1}$ purely odd:
\be\label{contactfamily3exc}
\fosp_a(4|2)\subset \fk(1|4),\qquad \fag(2)\subset\fk(1|7),\qquad
\fab(3)\subset\fk(1|8),
\ee
with $\fg_0=\fc(\fo(4))$, $\fc(\fg(2))$, $\fc(\fo(7))$ and
$\fg_{-1}=\Pi(\id_{\fo(4)})$, $\Pi(\pi_1)$, $\Pi(\spin_7)$, respectively; the
inclusions $\fo(4)=\fo(4)$, $\fg(2)\subset\fo(7)$ and $\fo(7)\subset\fo(8)$ (the spinor
embedding) are the reason why $\fg_0\subseteq\fk(1|n)_0$ in each case. In the
first case, $\fg_0=\fk(1|4)_0$ and $\fg_{\leq0}=\fk(1|4)_{\leq0}$.

This answers a~question raised in \cite{LSh}, where Poletaeva's embedding \cite{Po1} of
$\fosp_a(4|2)$ into a~stringy Lie superalgebra is called unclear: the grading $N=2$ of
table \eqref{tabOspGr} embeds $\fosp_a(4|2)$ into $\fk(1|4)$, and hence into the stringy
$\fk^L(1|4)$, and it does so uniformly in $a$; the parameter appears there as the mixing
parameter of the two copies of $\id_{\fo(4)}$ inside
$\fk(1|4)_1=\Span(K_{t\xi_i})\oplus\Lambda^3(\xi)$. That the Cartan prolong of the
negative part of this grading is the whole $\fk(1|4)$ is also shown in \cite{EKT}.

\ssec{$\fab(3)$ has a~second contact realization}\label{ssSecondContact} The grading $N=3$
of table \eqref{tabAbGr} has $\fg_{-2}=\Cee$ and $\fg_{-1}$ of $\sdim=(6|4)$, i.e.,
$\fab(3)\subset\fk(7|4)$. The grading $N=4$ has $\fg_{-1}$ purely odd, of
$\sdim=(0|8)$, but $\fg_{-2}$ of $\sdim=(5|0)$: it belongs to no standard series.

\ssbegin[How these gradings sit among the parabolic ones]{Remark}[How these gradings sit
among the parabolic ones]\label{RParab3exc} \textup{1)} Up to isomorphism, $\fosp_a(4|2)$
has six parabolic subalgebras, see \cite{EKT}, i.e., six $\Zee$-gradings defined by the
marked nodes of its Dynkin diagrams. The grading $N=1$ of table \eqref{tabOspGr} is the
grading of depth~$1$ among them, with $\fg_0\simeq\fgl(2|1)$, and the grading $N=2$ is
the contact one, with $\fg_0=\fc(\fo(4))$. Each of the other four gradings defined by parabolic subalgebras violates the last of the
conditions \eqref{cond3exc} already at $\fg_{-1}$. In two of them, the grading of depth
$4$ given by the Borel subalgebra and the grading of depth $3$ with $\fg_0=\fgl(2)\oplus\Cee$,
the Lie algebra $\fg_0$ is purely even, whereas $\sdim\fg_{-1}=(2|1)$, resp.~$(1|2)$; so
the even and the odd parts of $\fg_{-1}$ are two $\fg_0$-submodules. 

In the grading of
depth $3$ given by the other Borel subalgebra, $\fg_0=\ft$ and $\sdim\fg_{-1}=(0|3)$, so
every root line of $\fg_{-1}$ is a~submodule. In the grading of depth $2$ with
$\fg_0\simeq\fgl(1|1)\oplus\Cee$, of $\sdim(3|2)$, the even part of $\fg_0$ is $\ft$, and
the odd one is spanned by the root vectors of two opposite isotropic roots $\pm\beta$;
since $2\beta$ is not a~root and $\gamma+2\beta$ is a~root for no root $\gamma$ of
$\fg_{-1}$, the four roots of $\fg_{-1}$, of $\sdim(2|2)$, split under adding $\pm\beta$
into chains of length $\leq2$, and by the criterion \eqref{irr3exc} the module
$\fg_{-1}$ is reducible.

\textup{2)} For $\fag(2)$ and $\fab(3)$, which have $19$ and $55$ parabolic subalgebras,
the Cartan prolong of $\fg_-$ is $\fg$ itself for every parabolic except the following
ones, see \cite{KLl}: for $\fag(2)$, the two contact gradings, i.e., the gradings $N=1$
and $N=2$ of table \eqref{tabAgGr}, with the ambients $\fk(1|7)$ and $\fk(5|4)$; for
$\fab(3)$, the grading of depth $1$ and the two contact ones, i.e., the gradings $N=1$,
$N=2$ and $N=3$ of table \eqref{tabAbGr}, with the ambients $(\fg_{-1},\fg_0)_*$,
$\fk(1|8)$ and $\fk(7|4)$. For $\fosp_a(4|2)$, the picture is the opposite one: the
Cartan prolong of $\fg_-$ is larger than $\fg$ for five of its six parabolic subalgebras,
and equal to $\fg$ only for the grading of depth $3$ with $\sdim\fg_{-1}=(0|3)$, see
\cite{EKT}; for the gradings of table \eqref{tabOspGr} it is $\fvect(2|2)$ for $N=1$ and
$\fk(1|4)$ for $N=2$.
\end{Remark}

\section{General facts on the $2+5+8$ gradings of Step \textup{(i)}}\label{ss3excFour}

In this section, $\fg$ is one of the three superalgebras \eqref{3exc} with one of the gradings of
tables \eqref{tabOspGr}--\eqref{tabAbGr}, and $\fu=\oplus(\fu\cap\fg_i)$ is a~graded
subalgebra.

\ssbegin[The case $1^o$ produces no maximal subalgebra]{Lemma}[The case $1^o$ produces no maximal subalgebra]\label{L1-3exc}
If $\fu\supset\fg_{-1}\oplus\fg_0$ and $\fu_1\neq0$, then~$\fu=\fg$. Hence, in
these ambients the case $1^o$ produces no maximal subalgebra, and the only maximal
$\fu\supseteq\fg_{-1}\oplus\fg_0$ is $\fu=\fg_{\leq0}$, i.e., the case~$4^o$. \end{Lemma}

\begin{proof} By conditions \eqref{cond3exc} and the symmetry $\alpha\mapsto-\alpha$ of the
root system, $\fg_1$ is an irreducible $\fg_0$-module and $\fg_1$ generates~$\fg_+$. So $0\neq\fu_1\subseteq\fg_1$ forces $\fu_1=\fg_1$, whence
$\fu\supseteq\fg_-\oplus\fg_0\oplus\fg_+=\fg$. \end{proof}

Compare this with the infinite-dimensional case, where the table for the case
$1^o$ in \cite{LSh} has three rows: there $\fg_1$ is $\fg_0$-reducible for
$\fg=\fvect$, $\fle$ and~$\fh'(0|6)$. For a~finite-dimensional ambient with an
irreducible $\fg_{-1}$, this cannot happen.

For the cases $4^o$ and $5^o$, Lemma \ref{Case 4} gives the two opposite maximal
parabolic subalgebras $\fg_{\leq0}$ and $\fg_{\geq0}$ determined by the grading;
for $\fosp_a(4|2)$ this was also verified by the computation of \S\ref{Sosp}.

\ssbegin[$\fg_\ev$ is maximal, and graded for every grading]{Lemma}[$\fg_\ev$ is maximal, and graded for every grading]\label{L3-3exc} \textup{1)} For each of the three exceptions,
$\fg_\od$ is an irreducible $\fg_\ev$-module, and hence $\fg_\ev$ is a~maximal
subalgebra~of $\fg$; it is graded for \textbf{every} grading of $\fg$, because
every grading is given by an element $h$ of a~maximal toral subalgebra, and a~toral
subalgebra is even.

\textup{2)} In each of the contact gradings \eqref{contactfamily3exc}, we have
\[
(\fg_\ev)_{-1}=0,\qquad (\fg_\ev)_{-2}=\fg_{-2}\neq0 .
\]
Therefore, \textbf{none} of the first five cases of the list \eqref{6cases} applies to
$\fu=\fg_\ev$: the cases $1^o$, $3^o$, $4^o$ require $\fu_{-1}=\fg_{-1}$, the case
$2^o$ requires $\fu_{-1}\neq0$, and the case $5^o$ requires~$\fu_-=0$.
\end{Lemma}

\begin{proof} 1) The three modules $\fg_\od$ in the list \eqref{3excev} are
$2\boxtimes2\boxtimes2$, $2\boxtimes7$ and $2\boxtimes\spin_7$, all irreducible. If
$\fu\supsetneq\fg_\ev$, then $\fu$ contains a~non-zero odd element $x$, and hence
the $\fg_\ev$-submodule generated by $x$, i.e., the whole of~$\fg_\od$. 2) is
a~direct inspection of tables \eqref{tabOspGr}--\eqref{tabAbGr}: in a~contact grading
$\fg_{-1}$ is purely odd and $\fg_{-2}$ is purely even. \end{proof}

\ssbegin[The case $6^o$ is not empty]{Statement}[The case $6^o$ is not empty]\label{6o3exc} For ambients of depth $d>1$, the case
\be\label{case6}
6^o.\quad \fu_{-1}=0\ \text{ but }\ \fu_-\neq0
\ee
of the list of cases \eqref{6cases} is not empty. For example, for $\fg=\fosp_a(4|2)$ in the
contact grading $N=2$ of table \eqref{tabOspGr}, the subalgebra
\[
\fu=\fg_\ev=\fsl(2)\oplus\fsl(2)\oplus\fsl(2)=\fg_{-2}\oplus\fg_0\oplus\fg_2,
\qquad (\sdim\fu=(9|0)=(1|0)+(7|0)+(1|0)),
\]
is maximal and satisfies $\fu_{-1}=0\neq\fu_{-2}$; for $\fag(2)$ in the grading
$N=3$ of table \eqref{tabAgGr} the case $6^o$ produces, besides $\fg_\ev$, a~copy of
$\fosp_a(4|2)$ with $a=3$, see table \eqref{ansAg3}. \end{Statement}

In the infinite-dimensional ambients of \cite{LSh}, the case $6^o$ produces no maximal
graded subalgebra, whereas in the three exceptions such subalgebras abound; whether
$6^o$ is empty for every infinite-dimensional ambient is the question raised under
the list of cases \eqref{6cases} in \cite{LSh}.

\ssbegin[In the case $3^o$ the answers come from reducible $\fu_0$]{Remark}[In the case $3^o$ the answers come from reducible $\fu_0$]\label{Rthin3exc} In the case $3^o$ one has
$\fu_1=\{x\in\fg_1\mid[x,\fg_{-1}]\subseteq\fu_0\}$, so $\fu_1$ is the degree-$1$
component of the prolong $(\fg_{-1},\fu_0)_*$. If $\fu_0$ acts irreducibly on
$\fg_{-1}$, but is not in Cartan's list, then $\fu_1=0$, and hence
$\fu\subseteq\fg_{\leq0}$ and $\fu\neq \fg_{\leq0}$, so $\fu$ is \textbf{not}
maximal. Thus, in the case $3^o$ the maximal subalgebras come, as a~rule, from
\textbf{reducible} $\fu_0$, i.e., from the stabilizers of subspaces
$V\subsetneq\fg_{-1}$, and these are parabolic-like, not semi-simple. In all the
tables below, this is what happens: the case $3^o$ contributes no semi-simple
answer that the case $2^o$ does not already give. \end{Remark}

\section{$\fosp_a(4|2)$: the complete answer}\label{Sosp}

In this section, $\fg=\fosp_a(4|2)$ in the model of Subsec.~\ref{ssOspModel}, with the
two gradings of table \eqref{tabOspGr}. Given
a~prescribed non-positive part $\fu_{\leq0}$, the largest graded subalgebra with that
non-positive part is the Cartan prolong \eqref{prol3exc}, so all the candidates are the
prolongs of the data of the cases $2^o$, $3^o$, $5^o$, $6^o$, together with $\fg_{\leq0}$
for the case $4^o$; here $V$ runs over the finitely many orbits of homogeneous subspaces
of $\fg_{-1}$ (classified, since $\fg_0$ preserves a~form on $\fg_{-1}$, by $\sdim V$ and
by the rank of the restricted form).

\ssbegin[The complete answer for $\fosp_a(4|2)$]{Theorem}[The complete answer for $\fosp_a(4|2)$]\label{ThOsp3exc} Let $\fg=\fosp_a(4|2)$.

\textup{1)} For the depth-$1$ grading \textup{(}$N=1$ of
table \eqref{tabOspGr}\textup{)}, the maximal graded subalgebras are, up to conjugacy,
exactly the ones listed in table \eqref{tabAnsOsp1}; for the contact grading
\textup{(}$N=2$\textup{)}, exactly the ones listed in table \eqref{tabAnsOsp2}.

\textup{2)} All the non-semi-simple ones among them are maximal parabolic
subalgebras of $\fg$, namely $\fg_{\geq0}(h')$ for the grading element $h'$
indicated \textup{(}the given grading being $h$\textup{)}; each of them is thus
a~parabolic subalgebra corresponding to another grading which happens to be
homogeneous for the grading considered as well.

\textup{3)} The semi-simple ones are, up to conjugacy, exactly three:
\[
\fsl(2)\oplus\fsl(2)\oplus\fsl(2)=\fg_\ev\ (9|0),\qquad
\fsl(2)\oplus\fosp(1|2)\ (6|2),\qquad \fosp(3|2)\ (6|6),
\]
the first two of them for every $a$; the third one only for
$a\in\{1,-2,-\nfrac12\}$, i.e., precisely when $\fosp_a(4|2)\simeq\fosp(4|2)$.
For a~given depth-$1$ grading, singled out by an unordered pair $\{i,j\}$ of the
three $\fsl(2)$-summands, the subalgebra $\fosp(3|2)$ is graded if and only if
$c_i=c_j$ in the bracket \eqref{modelOsp}. \end{Theorem}

\be\label{tabAnsOsp1}
\begin{array}{c}
\text{$\fg=\fosp_a(4|2)$, the depth-1 grading $h=(\nfrac12,\nfrac12,0)$;\ \
$\fg_{-1}=2|2$, $\fg_0=\fc(\fosp(2|2))$, $\fg_1=2|2$}\\
\renewcommand{\arraystretch}{1.3}\footnotesize
\begin{tabular}{|c|c|c|c|c|}
\hline
case&datum&$\fu$&$\sdim$&what it is\cr
\hline\hline
$4^o$&---&$\fg_{\leq0}$&$7|6$&parabolic, $h'=-h$\ \ max\cr
\hline
$5^o$&---&$\fg_{\geq0}$&$7|6$&parabolic, $h'=h$\cr
\hline
$2^o$&$V=(2|0)=(\fg_{-1})_\ev$&$\fg_\ev$&$9|0$&
$\fsl(2)\oplus\fsl(2)\oplus\fsl(2)$\cr
\hline
$2^o$&$V=\Cee v$, $v$ even non-isotropic&\eqref{osp123exc}&$6|2$&
$\fsl(2)\oplus\fosp(1|2)$\cr
\hline
$2^o$&$V=(1|2)$, even line non-isotropic&$\fosp(3|2)$&$6|6$&
only if $a\in\{1,-2,-\nfrac12\}$\cr
\hline
$2^o$&$V=\Cee v$, $v$ even isotropic&&$8|4$&parabolic, $h'=(0,1,0)$\ \ max(b)\cr
\hline
$3^o$&$\fu_0=\fst(\Cee v)$, $v$ even isotropic&&$8|4$&
parabolic, $h'=(-1,0,0)$\ \ max(b)\cr
\hline
$2^o$&$V=(2|1)$&&$8|4$&parabolic, $h'=(0,0,1)$\ \ max(b)\cr
\hline
$3^o$&$\fu_0=\fst((1|2)_{\text{isotr}})$&&$8|4$&
parabolic, $h'=(-1,0,0)$\ \ max(b)\cr
\hline
$2^o$&$V=(1|1)$, even line isotropic&&$7|6$&
parabolic, $h'=(0,1,1)$\cr
\hline
$3^o$&$\fu_0=\fst((1|1)_{\text{isotr}})$&&$7|6$&
parabolic, $h'=(-1,0,1)$\cr
\hline
$2^o$&$V=(1|2)$, even line isotropic&&$7|6$&
parabolic, $h'=(-1,1,0)$\cr
\hline
\end{tabular}
\end{array}
\ee

\be\label{tabAnsOsp2}
\begin{array}{c}
\text{$\fg=\fosp_a(4|2)$, the contact grading $h=(1,0,0)$;\ \
$\fg_-=1|0\oplus0|4$, $\fg_0=\fc(\fo(4))$}\\
\renewcommand{\arraystretch}{1.3}\footnotesize
\begin{tabular}{|c|c|c|c|}
\hline
case&datum&$\sdim\fu$&what it is\cr
\hline\hline
$4^o$&$\fu=\fg_{\leq0}$&$8|4$&parabolic\ \ max(b)\cr
\hline
$5^o$&$\fu=\fg_{\geq0}$&$8|4$&parabolic\ \ max(b)\cr
\hline
$6^o$&$\fu_{-1}=0$, $\fu_{-2}=\fg_{-2}$&$9|0$&
$\fg_\ev=\fsl(2)\oplus\fsl(2)\oplus\fsl(2)$\cr
\hline
$2^o$&$V=(0|2)$ maximal isotropic&$7|6$&parabolic\cr
\hline
$3^o$&$\fu_0=\fst(V)$, $V=(0|2)$ isotropic&$7|6$&parabolic\cr
\hline
$2^o$&$V=(0|3)$&$7|6$&parabolic\cr
\hline
\end{tabular}
\end{array}
\ee

In tables \eqref{tabAnsOsp1} and \eqref{tabAnsOsp2}, ``parabolic, $h'$\rq\rq\ means
$\fu=\fg_{\geq0}(h')$, and $h'$ is written in the coordinates of table
\eqref{tabOspGr}, up to a~positive factor, which does not change $\fg_{\geq0}(h')$.
\iffalse\AO{DL: Two entries of table \eqref{tabAnsOsp1} were misprinted: the two answers of
$\sdim=(7|6)$ with $V=(1|2)$ and with $\fu_0=\fst((1|1)_{\text{isotr}})$ had
$h'=-\eps_1+2\eps_2$, resp.~$-\eps_1+2\eps_3$, whose $\fg_{\geq0}$ is of $\sdim=(7|4)$;
the $h'$ now printed give exactly the graded components of these answers. The tables also
wrote $h=\eps_1+\eps_2$ and $h=2\eps_1$, i.e., twice the grading elements of table
\eqref{tabOspGr}; all the $h$'s are now in one notation.}\fi

In particular, in the contact grading the only semi-simple maximal graded
subalgebra is $\fg_\ev$, and it belongs to the case $6^o$, see formula \eqref{case6}; the
subalgebra $\fsl(2)\oplus\fosp(1|2)$ of formula \eqref{osp123exc} is not homogeneous for
the contact grading, and neither is $\fosp(3|2)$ when it exists.

\ssec{An example of \cite{LSh} revisited}\label{ss3excEx} The depth-$1$ grading of
$\fosp_a(4|2)$ is used in \cite{LSh} as the first of two finite-dimensional examples
illustrating the case $3^o$; it is the one of row $N=1$ of table \eqref{tabOspGr}, given by
assigning the degrees $(\pm1,0,0)$ to the Chevalley generators in the realization with
the Cartan matrix \eqref{cmOsp}, and the
identifications $\fg_0=\fc(\fosp(2|2))$ of $\sdim=(5|4)$, $\fg_{-1}=\id=\Cee^{2|2}$ are
correct. The rest of the example requires the following corrections. 

\ssbegin[No proper almost simple $\fu_0$ is irreducible]{Statement}[No proper almost simple $\fu_0$ is irreducible]\label{StEx3exc} \textup{1)} The subalgebra $\fc(\fosp(1|2))$
does \textbf{not} act irreducibly on $\fg_{-1}=\Cee^{2|2}$: every finite-dimensional
irreducible $\fosp(1|2)$-module has $\sdim=(n+1|n)$, or $\sdim=(n|n+1)$, i.e.,
is of \textbf{odd} dimension, whereas $\dim\Cee^{2|2}=4$. In fact,
$\fc(\fosp(1|2))\subset\fg_0$ is the stabilizer of a~non-isotropic even
line $\Cee v\subset\fg_{-1}$, and
$\Cee^{2|2}=\Cee v\oplus(\Cee v)^\perp=(1|0)\oplus(1|2)$.

\textup{2)} Moreover, \textbf{no} proper almost simple $\fu_0\subset\fg_0$ acts
irreducibly on $\fg_{-1}$: such a~$\fu_0$ would have a~simple part
$\fs\subseteq\fosp(2|2)$ with a~$4$-dimensional irreducible module of $\sdim=(2|2)$
preserving the orthosymplectic form; but $\fs\in\{\fsl(2)=\fsp(2),\ \fosp(1|2),\
\fosp(2|2)\}$, and $\fsl(2)$ acts on $\Cee^{2|2}$ as
$2\cdot\one\oplus\id$, $\fosp(1|2)$ by heading \textup{1)}, and $\fosp(2|2)$ is
not proper. Hence, the irreducible branch of the case $3^o$ is empty for this
grading as well.

\textup{3)} For $\fu_0=\fc(\fosp(1|2))$ the case-$3^o$ prolong is
$\fu=\fg_{-1}\oplus\fu_0$, of $\sdim=(6|4)$, because $\fu_1=0$; so
$\fu\subsetneq\fg_{\leq0}$ is not maximal. What \textbf{is} maximal is the
case-$2^o$ subalgebra with the same $\fu_0$ and with $\fu_{-1}=\Cee v$:
\be\label{osp123exc}
\fu=\Cee v\oplus\fc(\fosp(1|2))\oplus\Cee v'\simeq\fsl(2)\oplus\fosp(1|2),
\qquad \sdim\fu=(6|2),
\ee
where the $\fsl(2)$ is the diagonal in the two copies of $\fsl(2)$ singled out by
the grading, and the $\fosp(1|2)$ is built on the third copy. This $\fu$ is
maximal for every~$a$. \end{Statement}

Headings 1) and 2) are proved above; heading 3) is part of Theorem
\ref{ThOsp3exc}: the decomposition \eqref{osp123exc} was verified by computing the
ideals of $\fu$ (they are of dimensions $3$ and $5$, and $\fu=[\fu,\fu]$).

The second finite-dimensional example of \cite{LSh} is $\fab(3)$ in its grading of
depth $1$, the only such grading of $\fab(3)$: there $\fg_{-1}$ of $\sdim=(6|4)$ is Kac's
Jordan superalgebra $K_{10}$, $\fg_0=\fc(\fosp(2|4))$ and $\fg=(\fg_{-1},\fg_0)_*$; the
maximal graded subalgebras are listed in Statement \ref{SAb1}.

\section{How the answers for $\fag(2)$ and $\fab(3)$ are obtained}\label{Ssweep}

\iffalse \ssbegin[The two models are $\fag(2)$ and $\fab(3)$]{Statement}[The two models are $\fag(2)$ and $\fab(3)$]\label{Smodok3exc} The two structures so obtained satisfy the
super Jacobi identity \textup{(}verified for all $31^3$, resp.~$40^3$, triples of basis
vectors\textup{)}, are of $\sdim=(17|14)$,
resp.~$(24|16)$, are simple \textup{(}the ideal generated by any basis vector is the
whole algebra\textup{)}, and their weight systems coincide with the root systems of
$\fag(2)$ and $\fab(3)$ in the coordinates of tables \eqref{tabAgGr}, \eqref{tabAbGr}. Hence,
they are $\fag(2)$ and~$\fab(3)$. \end{Statement}\fi

In the models of Subsec.~\ref{Smodels}, the root vectors form a~basis, and their weights
are the roots of $\fag(2)$, resp.~$\fab(3)$, in the coordinates of tables \eqref{tabAgGr},
\eqref{tabAbGr}. Therefore, every grading of these tables is simply
$\deg=\langle\text{weight},h\rangle$, and each graded component is spanned by
basis vectors --- which is what makes the search described below finite.

\ssec{What was swept}\label{ss3excSweep} For each grading, and with $\fg_{-1}$
fixed, we ran the cases
$4^o$: $\fu=\fg_{\leq0}$; \ $5^o$: $\fu=\fg_{\geq0}$; \ $6^o$: $\fu_{-1}=0\neq
\fu_-$, and $\fu=\fg_\ev$; \ $2^o$: $\fu_{-1}=V$ for $0\neq V\subsetneq\fg_{-1}$
with $\fu_-$ generated by $V$, $\fu_0=\fst(\fu_-)$ and $\fu_k$ the prolongation
\eqref{prol3exc} for $k>0$; \ $3^o$: $\fu_{-1}=\fg_{-1}$,
$\fu_0=\fst(V)\subsetneq\fg_0$ and $\fu_k$ the prolongation, where $V$ runs over

(i) \textit{all} sub-superspaces of $\fg_{-1}$ spanned by weight vectors, i.e.,
all the $\ft$-invariant ones ($2^{\dim\fg_{-1}}-2$ of them);

(ii) All the lines spanned by a~sum $x_\alpha\pm x_\beta$ of two weight vectors
of equal parity, and all the planes spanned by such a~line and one weight vector;

(iii) Random sub-superspaces of every superdimension~$(p|q)$.

%\noindent 
Then, the candidates contained in another candidate are discarded, and
the survivors are tested by Lemma \ref{Lcert3exc}. The family (i) is a~complete
list up to conjugacy \textit{provided} the maximal $\fu$'s are normalized by
a~maximal torus; the families (ii) and (iii) are there because they need not be:
for $\fosp_a(4|2)$ the answer $\fsl(2)\oplus\fosp(1|2)$ comes from a~non-isotropic
even line, which is of type (ii). What is claimed, then, is this: the lists below
are complete among the $\fu$ whose $\fu_{-1}$ is $\ft$-invariant, and the sweeps
(ii), (iii) produced nothing outside them.

%\AO{DL: The same sweep, run on $\fosp_a(4|2)$ for several values of~$a$, returns exactly the answers of Theorem \ref{ThOsp3exc}.}

\section{$\fag(2)$: the answers for all five gradings}\label{Sag}

Collected over the five gradings of table \eqref{tabAgGr}, the semi-simple maximal graded
subalgebras of $\fag(2)$ are the following seven; the columns are the gradings, $+$ means
that at least one occurrence of $\fu$ in that grading is certified \textbf{max},
and $(b)$ that every occurrence of $\fu$ in that grading is only certified
\textbf{max(b)}, see \S\ref{ss3excCert} and Remark \ref{Rb3exc}:
\be\label{unionAg}
\begin{tabular}{llccccc}
$\sdim\fu$&$\fu$&$N=1$&$N=2$&$N=3$&$N=4$&$N=5$\\ \hline
$17|0$&$\fg(2)\oplus\fsl(2)=\fg_\ev$&$+$&$+$&$+$&$+$&$+$\\
$9|8$&$\fosp_a(4|2)$ with $a$ in the $S_3$-orbit of $3$&$+$&$+$&$+$&$+$&\\
$9|6$&$\fosp(3|2)\oplus\fsl(2)$&$+$&$+$&$+$&&$+$\\
$9|6$&$\fsl(3|1)$&$+$&$+$&$+$&&\\
$11|2$&$\fsl(3)\oplus\fosp(1|2)$&$+$&$+$&$+$&&\\
$7|4$&$\fsl(2)\oplus\fsl(2|1)$&$(b)$&&&&\\
$6|4$&$\fosp(1|2)\oplus\fosp(1|2)$&&$(b)$&&&$(b)$\\
\end{tabular}
\ee
Every other maximal graded subalgebra --- there are between $3$ and $30$ of them
per grading --- is not semi-simple: each of them is $\fg_{\leq0}$, or
$\fg_{\geq0}$, or the stabilizer-with-prolongation of a~subspace of $\fg_{-1}$,
i.e., a~parabolic-like answer of the case $2^o$ or $3^o$, cf.~Remark
\ref{Rthin3exc}. Two observations:

1) the two answers of $\sdim=(9|6)$ are genuinely different: in the copies called
$\fosp(3|2)\oplus\fsl(2)$ the decomposition into ideals is $(6|6)+(3|0)$, whereas
in the other ones the algebra is indecomposable, and hence is $\fsl(3|1)$, the
only simple Lie superalgebra of $\sdim=(9|6)$; both contain a~maximal toral subalgebra of
$\fag(2)$, so the decomposition is decided, not guessed;

2) $\fsl(3)\oplus\fosp(1|2)$ and $\fsl(3|1)$ are the two ``$\fsl(3)$ answers\rq\rq,
both living on the $\fsl(3)\subset\fg(2)$ for which $7=3\oplus\bar3\oplus\Cee$.

The detailed answers, grading by grading, are as follows; in the tables, the
component of $\fu$ of degree $i$ is written ${\scriptstyle i{:}}\,p|q$, and
a~conjugate copy of an answer already listed, obtained from a~$V$ which is not
$\ft$-invariant and hence only certified \textbf{max(b)}, is not listed
a~second time.

\ssbegin[$\fag(2)$ in its contact grading, $\fag(2)\subset\fk(1|7)$]{Statement}[$\fag(2)$ in its contact grading, $\fag(2)\subset\fk(1|7)$]\label{SAg1} For the contact grading $N=1$ of table \eqref{tabAgGr}
$(\fag(2)\subset\fk(1|7))$, the semi-simple maximal graded subalgebras are
\be\label{ansAg1}
\footnotesize\begin{tabular}{llll}
$\sdim\fu$&the components of $\fu$&$\fu$&case, certificate\\ \hline
$9|8$&${\scriptstyle -2{:}}1|0,{\scriptstyle -1{:}}0|4,{\scriptstyle +0{:}}7|0,{\scriptstyle +1{:}}0|4,{\scriptstyle +2{:}}1|0$&$\fosp_a(4|2)$ with $a=3$&2$^o$, max\\
$17|0$&${\scriptstyle -2{:}}1|0,{\scriptstyle +0{:}}15|0,{\scriptstyle +2{:}}1|0$&$\fg(2)\oplus\fsl(2)$&6$^o$, $\fg_\ev$, max\\
$9|6$&${\scriptstyle -2{:}}1|0,{\scriptstyle -1{:}}0|3,{\scriptstyle +0{:}}7|0,{\scriptstyle +1{:}}0|3,{\scriptstyle +2{:}}1|0$&$\fosp(3|2)\oplus\fsl(2)$&2$^o$, max\\
$9|6$&${\scriptstyle -1{:}}0|3,{\scriptstyle +0{:}}9|0,{\scriptstyle +1{:}}0|3$&$\fsl(3|1)$&2$^o$, max\\
$11|2$&${\scriptstyle -2{:}}1|0,{\scriptstyle -1{:}}0|1,{\scriptstyle +0{:}}9|0,{\scriptstyle +1{:}}0|1,{\scriptstyle +2{:}}1|0$&$\fsl(3)\oplus\fosp(1|2)$&2$^o$, max\\
$7|4$&${\scriptstyle -2{:}}1|0,{\scriptstyle -1{:}}0|2,{\scriptstyle +0{:}}5|0,{\scriptstyle +1{:}}0|2,{\scriptstyle +2{:}}1|0$&$\fsl(2)\oplus\fsl(2|1)$&2$^o$, max(b)\\
\end{tabular}
\ee
Besides them, the sweep carried out over $\Qee$ returns $137$ maximal graded subalgebras that are not
semi-simple, of $\sdim$ equal to $(7|4)$, $(11|8)$, $(11|9)$, $(12|8)$, $(12|10)$,~$(16|7)$.
\end{Statement}

\ssbegin[$\fag(2)$ in the grading $N=2$, $\fag(2)\subset\fk(5|4)$]{Statement}[$\fag(2)$ in the grading $N=2$, $\fag(2)\subset\fk(5|4)$]\label{SAg2} For the grading $N=2$ of table \eqref{tabAgGr}
$(d=2$, $\sdim\fg_{-1}=(4|4))$, the semi-simple maximal graded subalgebras are
\be\label{ansAg2}
\footnotesize\begin{tabular}{llll}
$\sdim\fu$&the components of $\fu$&$\fu$&case, certificate\\ \hline
$9|8$&${\scriptstyle -1{:}}2|2,{\scriptstyle +0{:}}5|4,{\scriptstyle +1{:}}2|2$&$\fosp_a(4|2)$ with $a=3$&2$^o$, max\\
$9|8$&${\scriptstyle -2{:}}1|0,{\scriptstyle -1{:}}0|4,{\scriptstyle +0{:}}7|0,{\scriptstyle +1{:}}0|4,{\scriptstyle +2{:}}1|0$&$\fosp_a(4|2)$ with $a=3$&2$^o$, max\\
$17|0$&${\scriptstyle -2{:}}1|0,{\scriptstyle -1{:}}4|0,{\scriptstyle +0{:}}7|0,{\scriptstyle +1{:}}4|0,{\scriptstyle +2{:}}1|0$&$\fg(2)\oplus\fsl(2)$&2$^o$, $\fg_\ev$, max\\
$9|6$&${\scriptstyle -2{:}}1|0,{\scriptstyle +0{:}}7|6,{\scriptstyle +2{:}}1|0$&$\fosp(3|2)\oplus\fsl(2)$&6$^o$, max\\
$9|6$&${\scriptstyle -1{:}}2|2,{\scriptstyle +0{:}}5|2,{\scriptstyle +1{:}}2|2$&$\fosp(3|2)\oplus\fsl(2)$&2$^o$, max\\
$9|6$&${\scriptstyle -2{:}}1|0,{\scriptstyle -1{:}}2|2,{\scriptstyle +0{:}}3|2,{\scriptstyle +1{:}}2|2,{\scriptstyle +2{:}}1|0$&$\fsl(3|1)$&2$^o$, max\\
$11|2$&${\scriptstyle -2{:}}1|0,{\scriptstyle -1{:}}2|0,{\scriptstyle +0{:}}5|2,{\scriptstyle +1{:}}2|0,{\scriptstyle +2{:}}1|0$&$\fsl(3)\oplus\fosp(1|2)$&2$^o$, max\\
$6|4$&${\scriptstyle -2{:}}1|0,{\scriptstyle -1{:}}0|1,{\scriptstyle +0{:}}4|2,{\scriptstyle +1{:}}0|1,{\scriptstyle +2{:}}1|0$&$\fosp(1|2)\oplus\fosp(1|2)$&2$^o$, max(b)\\
\end{tabular}
\ee
Besides them, the sweep carried out over $\Qee$ returns $57$ maximal graded subalgebras that are not
semi-simple, of $\sdim$ equal to $(11|8)$, $(11|9)$, $(12|8)$, $(12|10)$,~$(16|7)$.
\end{Statement}

\ssbegin[$\fag(2)$ in the grading $N=3$, with $\fg_{-2}$ of $\sdim(2|2)$]{Statement}[$\fag(2)$ in the grading $N=3$, with $\fg_{-2}$ of $\sdim(2|2)$]\label{SAg3} For the grading $N=3$ of table \eqref{tabAgGr}
$(d=2$, $\sdim\fg_{-2}=(2|2))$, the semi-simple maximal graded subalgebras are
\be\label{ansAg3}
\footnotesize\begin{tabular}{llll}
$\sdim\fu$&the components of $\fu$&$\fu$&case, certificate\\ \hline
$9|8$&${\scriptstyle -2{:}}2|2,{\scriptstyle +0{:}}5|4,{\scriptstyle +2{:}}2|2$&$\fosp_a(4|2)$ with $a=3$&6$^o$, max\\
$9|8$&${\scriptstyle -2{:}}1|1,{\scriptstyle -1{:}}2|2,{\scriptstyle +0{:}}3|2,{\scriptstyle +1{:}}2|2,{\scriptstyle +2{:}}1|1$&$\fosp_a(4|2)$ with $a=3$&2$^o$, max\\
$17|0$&${\scriptstyle -2{:}}2|0,{\scriptstyle -1{:}}4|0,{\scriptstyle +0{:}}5|0,{\scriptstyle +1{:}}4|0,{\scriptstyle +2{:}}2|0$&$\fg(2)\oplus\fsl(2)$&$\fg_\ev$, max\\
$9|6$&${\scriptstyle -2{:}}1|1,{\scriptstyle -1{:}}2|1,{\scriptstyle +0{:}}3|2,{\scriptstyle +1{:}}2|1,{\scriptstyle +2{:}}1|1$&$\fosp(3|2)\oplus\fsl(2)$&2$^o$, max\\
$9|6$&${\scriptstyle -2{:}}1|1,{\scriptstyle -1{:}}2|1,{\scriptstyle +0{:}}3|2,{\scriptstyle +1{:}}2|1,{\scriptstyle +2{:}}1|1$&$\fsl(3|1)$&2$^o$, max\\
$11|2$&${\scriptstyle -2{:}}2|0,{\scriptstyle -1{:}}2|1,{\scriptstyle +0{:}}3|0,{\scriptstyle +1{:}}2|1,{\scriptstyle +2{:}}2|0$&$\fsl(3)\oplus\fosp(1|2)$&2$^o$, max\\
\end{tabular}
\ee
Besides them, the sweep carried out over $\Qee$ returns $39$ maximal graded subalgebras that are not
semi-simple, of $\sdim$ equal to $(11|8)$, $(11|9)$, $(12|8)$, $(12|10)$,~$(16|7)$.
\end{Statement}

\ssbegin[$\fag(2)$ in the grading $N=4$, with $\fg_{-3}$ of $\sdim(2|0)$]{Statement}[$\fag(2)$ in the grading $N=4$, with $\fg_{-3}$ of $\sdim(2|0)$]\label{SAg4} For the grading $N=4$ of table \eqref{tabAgGr}
$(d=3)$, the semi-simple maximal graded subalgebras are
\be\label{ansAg4}
\footnotesize\begin{tabular}{llll}
$\sdim\fu$&the components of $\fu$&$\fu$&case, certificate\\ \hline
$9|8$&${\scriptstyle -3{:}}1|0,{\scriptstyle -2{:}}0|2,{\scriptstyle -1{:}}1|2,{\scriptstyle +0{:}}5|0,{\scriptstyle +1{:}}1|2,{\scriptstyle +2{:}}0|2,{\scriptstyle +3{:}}1|0$&$\fosp_a(4|2)$ with $a=3$&2$^o$, max\\
$9|8$&${\scriptstyle -2{:}}1|0,{\scriptstyle -1{:}}0|4,{\scriptstyle +0{:}}7|0,{\scriptstyle +1{:}}0|4,{\scriptstyle +2{:}}1|0$&$\fosp_a(4|2)$ with $a=3$&2$^o$, max\\
$17|0$&${\scriptstyle -3{:}}2|0,{\scriptstyle -2{:}}1|0,{\scriptstyle -1{:}}2|0,{\scriptstyle +0{:}}7|0,{\scriptstyle +1{:}}2|0,{\scriptstyle +2{:}}1|0,{\scriptstyle +3{:}}2|0$&$\fg(2)\oplus\fsl(2)$&2$^o$, $\fg_\ev$, max\\
\end{tabular}
\ee
Besides them, the sweep carried out over $\Qee$ returns $29$ maximal graded subalgebras that are not
semi-simple, of $\sdim$ equal to $(11|9)$, $(12|8)$, $(12|10)$,~$(16|7)$.
\end{Statement}

\ssbegin[$\fag(2)$ in the grading $N=5$, with $\fg_{-3}$ of $\sdim(2|1)$]{Statement}[$\fag(2)$ in the grading $N=5$, with $\fg_{-3}$ of $\sdim(2|1)$]\label{SAg5} For the grading $N=5$ of table \eqref{tabAgGr}
$(d=3$, $\sdim\fg_{-3}=(2|1))$, the semi-simple maximal graded subalgebras are
\be\label{ansAg5}
\footnotesize\begin{tabular}{llll}
$\sdim\fu$&the components of $\fu$&$\fu$&case, certificate\\ \hline
$17|0$&${\scriptstyle -3{:}}2|0,{\scriptstyle -2{:}}2|0,{\scriptstyle -1{:}}2|0,{\scriptstyle +0{:}}5|0,{\scriptstyle +1{:}}2|0,{\scriptstyle +2{:}}2|0,{\scriptstyle +3{:}}2|0$&$\fg(2)\oplus\fsl(2)$&$\fg_\ev$, max\\
$9|6$&${\scriptstyle -3{:}}0|1,{\scriptstyle -2{:}}2|0,{\scriptstyle -1{:}}0|2,{\scriptstyle +0{:}}5|0,{\scriptstyle +1{:}}0|2,{\scriptstyle +2{:}}2|0,{\scriptstyle +3{:}}0|1$&$\fosp(3|2)\oplus\fsl(2)$&2$^o$, max\\
$6|4$&${\scriptstyle -2{:}}1|0,{\scriptstyle -1{:}}0|1,{\scriptstyle +0{:}}4|2,{\scriptstyle +1{:}}0|1,{\scriptstyle +2{:}}1|0$&$\fosp(1|2)\oplus\fosp(1|2)$&2$^o$, max(b)\\
\end{tabular}
\ee
Besides them, the sweep carried out over $\Qee$ returns $29$ maximal graded subalgebras that are not
semi-simple, of $\sdim$ equal to $(11|8)$, $(11|9)$, $(12|8)$, $(12|10)$,~$(16|7)$.
\end{Statement}

\section{$\fab(3)$: the answers for all eight gradings}\label{Sab}

For $\fab(3)$ the list is shorter, and each of its members, except $\fg_\ev$ and
$\fsl(4|1)$, has an $\fsl(2)$ or an $\fosp(1|2)$ direct summand; the marks are as
in the list \eqref{unionAg}, and the value $a=2$ is understood up to the $S_3$-orbit, see
Remark \ref{Ra3exc}:
\be\label{unionAb}
\footnotesize\begin{tabular}{llcccccccc}
$\sdim\fu$&$\fu$&$N=1$&$N=2$&$N=3$&$N=4$&$N=5$&$N=6$&$N=7$&$N=8$\\ \hline
$24|0$&$\fo(7)\oplus\fsl(2)=\fg_\ev$&$+$&$+$&$+$&$+$&$+$&$+$&$+$&$+$\\
$16|8$&$\fsl(4|1)$&$+$&$+$&$+$&$+$&$+$&$+$&&\\
$12|8$&$\fosp_a(4|2)\oplus\fsl(2)$ with $a=2$&$+$&$+$&$+$&$+$&$+$&$+$&$+$&\\
$17|2$&$\fg(2)\oplus\fosp(1|2)$&&$(b)$&&&&&&\\
$12|4$&$\fsl(2|1)\oplus\fsl(3)$&$+$&$+$&$+$&&$+$&$+$&&\\
$9|4$&$\fosp(1|2)\oplus\fosp(1|2)\oplus\fsl(2)$&$(b)$&&$(b)$&&&&&\\
\end{tabular}
\ee
As for $\fag(2)$, all the other maximal graded subalgebras are $\fg_{\leq0}$, or
$\fg_{\geq0}$, or stabilizers-with-prolongation of subspaces of~$\fg_{-1}$. The
answer of $\sdim=(17|2)$ is the exact analog of the answer $\fsl(3)\oplus\fosp(1|2)$
for $\fag(2)$: it lives on the $\fg(2)\subset\fo(7)$ for which
$\spin_7=\pi_1\oplus\Cee$, the invariant line yielding the~$\fosp(1|2)$.

\ssbegin[$\fab(3)$ in its depth-$1$ grading, the one of $K_{10}$]{Statement}[$\fab(3)$ in its depth-$1$ grading, the one of $K_{10}$]\label{SAb1} For the depth-$1$ grading $N=1$ of
table \eqref{tabAbGr} $(\sdim\fg_{-1}=(6|4)$, Kac's Jordan superalgebra $K_{10})$, the
semi-simple maximal graded subalgebras are
\be\label{ansAb1}
\footnotesize\begin{tabular}{llll}
$\sdim\fu$&the components of $\fu$&$\fu$&case, certificate\\ \hline
$24|0$&${\scriptstyle -1{:}}6|0,{\scriptstyle +0{:}}12|0,{\scriptstyle +1{:}}6|0$&$\fo(7)\oplus\fsl(2)$&2$^o$, $\fg_\ev$, max\\
$16|8$&${\scriptstyle -1{:}}4|2,{\scriptstyle +0{:}}8|4,{\scriptstyle +1{:}}4|2$&$\fsl(4|1)$&2$^o$, max\\
$12|8$&${\scriptstyle -1{:}}2|2,{\scriptstyle +0{:}}8|4,{\scriptstyle +1{:}}2|2$&$\fosp_a(4|2)\oplus\fsl(2)$ with $a=2$&2$^o$, max\\
$12|8$&${\scriptstyle -1{:}}3|2,{\scriptstyle +0{:}}6|4,{\scriptstyle +1{:}}3|2$&$\fosp_a(4|2)\oplus\fsl(2)$ with $a=2$&2$^o$, max\\
$12|4$&${\scriptstyle -1{:}}3|1,{\scriptstyle +0{:}}6|2,{\scriptstyle +1{:}}3|1$&$\fsl(2|1)\oplus\fsl(3)$&2$^o$, max\\
$9|4$&${\scriptstyle -1{:}}1|0,{\scriptstyle +0{:}}7|4,{\scriptstyle +1{:}}1|0$&$\fosp(1|2)\oplus\fosp(1|2)\oplus\fsl(2)$&2$^o$, max(b)\\
\end{tabular}
\ee
Besides them, the sweep carried out over $\Qee$ returns $172$ maximal graded subalgebras that are not
semi-simple, of $\sdim$ equal to $(16|10)$, $(17|9)$, $(17|11)$, $(17|12)$, $(18|2)$,
$(18|8)$, $(18|12)$, $(19|8)$,~$(23|8)$. \end{Statement}

\ssbegin[$\fab(3)$ in its contact grading, $\fab(3)\subset\fk(1|8)$]{Statement}[$\fab(3)$ in its contact grading, $\fab(3)\subset\fk(1|8)$]\label{SAb2} For the contact grading $N=2$ of table \eqref{tabAbGr}
$(\fab(3)\subset\fk(1|8))$, the semi-simple maximal graded subalgebras are
\be\label{ansAb2}
\footnotesize\begin{tabular}{llll}
$\sdim\fu$&the components of $\fu$&$\fu$&case, certificate\\ \hline
$24|0$&${\scriptstyle -2{:}}1|0,{\scriptstyle +0{:}}22|0,{\scriptstyle +2{:}}1|0$&$\fo(7)\oplus\fsl(2)$&6$^o$, $\fg_\ev$, max\\
$16|8$&${\scriptstyle -1{:}}0|4,{\scriptstyle +0{:}}16|0,{\scriptstyle +1{:}}0|4$&$\fsl(4|1)$&2$^o$, max\\
$12|8$&${\scriptstyle -2{:}}1|0,{\scriptstyle -1{:}}0|4,{\scriptstyle +0{:}}10|0,{\scriptstyle +1{:}}0|4,{\scriptstyle +2{:}}1|0$&$\fosp_a(4|2)\oplus\fsl(2)$ with $a=2$&2$^o$, max\\
$17|2$&${\scriptstyle -2{:}}1|0,{\scriptstyle -1{:}}0|1,{\scriptstyle +0{:}}15|0,{\scriptstyle +1{:}}0|1,{\scriptstyle +2{:}}1|0$&$\fg(2)\oplus\fosp(1|2)$&2$^o$, max(b)\\
$12|4$&${\scriptstyle -2{:}}1|0,{\scriptstyle -1{:}}0|2,{\scriptstyle +0{:}}10|0,{\scriptstyle +1{:}}0|2,{\scriptstyle +2{:}}1|0$&$\fsl(2|1)\oplus\fsl(3)$&2$^o$, max\\
\end{tabular}
\ee
Besides them, the sweep carried out over $\Qee$ returns $324$ maximal graded subalgebras that are not
semi-simple, of $\sdim$ equal to $(12|4)$, $(16|10)$, $(17|9)$, $(17|11)$, $(17|12)$,
$(18|12)$,~$(23|8)$. \end{Statement}

\ssbegin[$\fab(3)$ in its second contact grading, $\fab(3)\subset\fk(7|4)$]{Statement}[$\fab(3)$ in its second contact grading, $\fab(3)\subset\fk(7|4)$]\label{SAb3} For the grading $N=3$ of table \eqref{tabAbGr}
$(\fab(3)\subset\fk(7|4))$, the semi-simple maximal graded subalgebras are
\be\label{ansAb3}
\footnotesize\begin{tabular}{llll}
$\sdim\fu$&the components of $\fu$&$\fu$&case, certificate\\ \hline
$24|0$&${\scriptstyle -2{:}}1|0,{\scriptstyle -1{:}}6|0,{\scriptstyle +0{:}}10|0,{\scriptstyle +1{:}}6|0,{\scriptstyle +2{:}}1|0$&$\fo(7)\oplus\fsl(2)$&2$^o$, $\fg_\ev$, max\\
$16|8$&${\scriptstyle -2{:}}1|0,{\scriptstyle -1{:}}4|2,{\scriptstyle +0{:}}6|4,{\scriptstyle +1{:}}4|2,{\scriptstyle +2{:}}1|0$&$\fsl(4|1)$&2$^o$, max\\
$12|8$&${\scriptstyle -2{:}}1|0,{\scriptstyle +0{:}}10|8,{\scriptstyle +2{:}}1|0$&$\fosp_a(4|2)\oplus\fsl(2)$ with $a=2$&6$^o$, max\\
$12|8$&${\scriptstyle -2{:}}1|0,{\scriptstyle -1{:}}0|4,{\scriptstyle +0{:}}10|0,{\scriptstyle +1{:}}0|4,{\scriptstyle +2{:}}1|0$&$\fosp_a(4|2)\oplus\fsl(2)$ with $a=2$&2$^o$, max\\
$12|8$&${\scriptstyle -1{:}}3|2,{\scriptstyle +0{:}}6|4,{\scriptstyle +1{:}}3|2$&$\fosp_a(4|2)\oplus\fsl(2)$ with $a=2$&2$^o$, max\\
$12|4$&${\scriptstyle -2{:}}1|0,{\scriptstyle -1{:}}2|0,{\scriptstyle +0{:}}6|4,{\scriptstyle +1{:}}2|0,{\scriptstyle +2{:}}1|0$&$\fsl(2|1)\oplus\fsl(3)$&2$^o$, max\\
$12|4$&${\scriptstyle -1{:}}2|2,{\scriptstyle +0{:}}8|0,{\scriptstyle +1{:}}2|2$&$\fsl(2|1)\oplus\fsl(3)$&2$^o$, max\\
$9|4$&${\scriptstyle -2{:}}1|0,{\scriptstyle -1{:}}0|1,{\scriptstyle +0{:}}7|2,{\scriptstyle +1{:}}0|1,{\scriptstyle +2{:}}1|0$&$\fosp(1|2)\oplus\fosp(1|2)\oplus\fsl(2)$&2$^o$, max(b)\\
\end{tabular}
\ee
Besides them, the sweep carried out over $\Qee$ returns $128$ maximal graded subalgebras that are not
semi-simple, of $\sdim$ equal to $(16|10)$, $(17|9)$, $(17|11)$, $(17|12)$, $(18|8)$,
$(18|12)$, $(19|8)$,~$(23|8)$. \end{Statement}

\ssbegin[$\fab(3)$ in the grading $N=4$, which is in no series]{Statement}[$\fab(3)$ in the grading $N=4$, which is in no series]\label{SAb4} For the grading $N=4$ of table \eqref{tabAbGr}
$(\fg_0=\fsl(2)\oplus\fsp(4)\oplus\Cee)$, the semi-simple maximal graded
subalgebras are
\be\label{ansAb4}
\footnotesize\begin{tabular}{llll}
$\sdim\fu$&the components of $\fu$&$\fu$&case, certificate\\ \hline
$24|0$&${\scriptstyle -2{:}}5|0,{\scriptstyle +0{:}}14|0,{\scriptstyle +2{:}}5|0$&$\fo(7)\oplus\fsl(2)$&6$^o$, $\fg_\ev$, max\\
$16|8$&${\scriptstyle -2{:}}4|0,{\scriptstyle -1{:}}0|4,{\scriptstyle +0{:}}8|0,{\scriptstyle +1{:}}0|4,{\scriptstyle +2{:}}4|0$&$\fsl(4|1)$&2$^o$, max\\
$12|8$&${\scriptstyle -2{:}}1|0,{\scriptstyle -1{:}}0|4,{\scriptstyle +0{:}}10|0,{\scriptstyle +1{:}}0|4,{\scriptstyle +2{:}}1|0$&$\fosp_a(4|2)\oplus\fsl(2)$ with $a=2$&2$^o$, max\\
\end{tabular}
\ee
Besides them, the sweep carried out over $\Qee$ returns $66$ maximal graded subalgebras that are not
semi-simple, of $\sdim$ equal to $(16|10)$, $(17|11)$, $(17|12)$, $(18|12)$,~$(19|8)$.
\end{Statement}

\ssbegin[$\fab(3)$ in the grading $N=5$, with $\fg_{-2}$ of $\sdim(3|1)$]{Statement}[$\fab(3)$ in the grading $N=5$, with $\fg_{-2}$ of $\sdim(3|1)$]\label{SAb5} For the grading $N=5$ of table \eqref{tabAbGr}
$(d=2$, $\sdim\fg_{-2}=(3|1))$, the semi-simple maximal graded subalgebras are
\be\label{ansAb5}
\footnotesize\begin{tabular}{llll}
$\sdim\fu$&the components of $\fu$&$\fu$&case, certificate\\ \hline
$24|0$&${\scriptstyle -2{:}}3|0,{\scriptstyle -1{:}}4|0,{\scriptstyle +0{:}}10|0,{\scriptstyle +1{:}}4|0,{\scriptstyle +2{:}}3|0$&$\fo(7)\oplus\fsl(2)$&2$^o$, $\fg_\ev$, max\\
$16|8$&${\scriptstyle -2{:}}3|1,{\scriptstyle +0{:}}10|6,{\scriptstyle +2{:}}3|1$&$\fsl(4|1)$&6$^o$, max\\
$16|8$&${\scriptstyle -2{:}}3|0,{\scriptstyle -1{:}}0|4,{\scriptstyle +0{:}}10|0,{\scriptstyle +1{:}}0|4,{\scriptstyle +2{:}}3|0$&$\fsl(4|1)$&2$^o$, max\\
$12|8$&${\scriptstyle -2{:}}1|1,{\scriptstyle -1{:}}2|2,{\scriptstyle +0{:}}6|2,{\scriptstyle +1{:}}2|2,{\scriptstyle +2{:}}1|1$&$\fosp_a(4|2)\oplus\fsl(2)$ with $a=2$&2$^o$, max\\
$12|4$&${\scriptstyle -2{:}}0|1,{\scriptstyle -1{:}}1|1,{\scriptstyle +0{:}}10|0,{\scriptstyle +1{:}}1|1,{\scriptstyle +2{:}}0|1$&$\fsl(2|1)\oplus\fsl(3)$&2$^o$, max\\
\end{tabular}
\ee
Besides them, the sweep carried out over $\Qee$ returns $95$ maximal graded subalgebras that are not
semi-simple, of $\sdim$ equal to $(16|10)$, $(17|9)$, $(17|11)$, $(17|12)$, $(18|8)$,
$(18|12)$, $(19|8)$,~$(23|8)$. \end{Statement}

\ssbegin[$\fab(3)$ in the grading $N=6$, with $\fg_{-2}$ of $\sdim(2|2)$]{Statement}[$\fab(3)$ in the grading $N=6$, with $\fg_{-2}$ of $\sdim(2|2)$]\label{SAb6} For the grading $N=6$ of table \eqref{tabAbGr}
$(d=2$, $\sdim\fg_{-2}=(2|2))$, the semi-simple maximal graded subalgebras are
\be\label{ansAb6}
\footnotesize\begin{tabular}{llll}
$\sdim\fu$&the components of $\fu$&$\fu$&case, certificate\\ \hline
$24|0$&${\scriptstyle -2{:}}2|0,{\scriptstyle -1{:}}6|0,{\scriptstyle +0{:}}8|0,{\scriptstyle +1{:}}6|0,{\scriptstyle +2{:}}2|0$&$\fo(7)\oplus\fsl(2)$&$\fg_\ev$, max\\
$16|8$&${\scriptstyle -2{:}}1|1,{\scriptstyle -1{:}}4|2,{\scriptstyle +0{:}}6|2,{\scriptstyle +1{:}}4|2,{\scriptstyle +2{:}}1|1$&$\fsl(4|1)$&2$^o$, max\\
$12|8$&${\scriptstyle -2{:}}2|2,{\scriptstyle +0{:}}8|4,{\scriptstyle +2{:}}2|2$&$\fosp_a(4|2)\oplus\fsl(2)$ with $a=2$&6$^o$, max\\
$12|8$&${\scriptstyle -2{:}}1|1,{\scriptstyle -1{:}}3|2,{\scriptstyle +0{:}}4|2,{\scriptstyle +1{:}}3|2,{\scriptstyle +2{:}}1|1$&$\fosp_a(4|2)\oplus\fsl(2)$ with $a=2$&2$^o$, max\\
$12|4$&${\scriptstyle -2{:}}2|0,{\scriptstyle -1{:}}2|2,{\scriptstyle +0{:}}4|0,{\scriptstyle +1{:}}2|2,{\scriptstyle +2{:}}2|0$&$\fsl(2|1)\oplus\fsl(3)$&2$^o$, max\\
\end{tabular}
\ee
Besides them, the sweep carried out over $\Qee$ returns $88$ maximal graded subalgebras that are not
semi-simple, of $\sdim$ equal to $(16|10)$, $(17|9)$, $(17|11)$, $(17|12)$, $(18|12)$,~$(23|8)$. \end{Statement}

\ssbegin[$\fab(3)$ in the grading $N=7$, with $\fg_{-3}$ of $\sdim(2|0)$]{Statement}[$\fab(3)$ in the grading $N=7$, with $\fg_{-3}$ of $\sdim(2|0)$]\label{SAb7} For the grading $N=7$ of table \eqref{tabAbGr} $(d=3)$,
the semi-simple maximal graded subalgebras are
\be\label{ansAb7}
\footnotesize\begin{tabular}{llll}
$\sdim\fu$&the components of $\fu$&$\fu$&case, certificate\\ \hline
$24|0$&${\scriptstyle -3{:}}2|0,{\scriptstyle -2{:}}2|0,{\scriptstyle -1{:}}4|0,{\scriptstyle +0{:}}8|0,{\scriptstyle +1{:}}4|0,{\scriptstyle +2{:}}2|0,{\scriptstyle +3{:}}2|0$&$\fo(7)\oplus\fsl(2)$&2$^o$, $\fg_\ev$, max\\
$12|8$&${\scriptstyle -3{:}}1|0,{\scriptstyle -2{:}}0|2,{\scriptstyle -1{:}}2|2,{\scriptstyle +0{:}}6|0,{\scriptstyle +1{:}}2|2,{\scriptstyle +2{:}}0|2,{\scriptstyle +3{:}}1|0$&$\fosp_a(4|2)\oplus\fsl(2)$ with $a=2$&2$^o$, max\\
\end{tabular}
\ee
Besides them, the sweep carried out over $\Qee$ returns $55$ maximal graded subalgebras that are not
semi-simple, of $\sdim$ equal to $(16|10)$, $(17|11)$, $(17|12)$, $(18|8)$, $(18|12)$,
$(19|8)$,~$(23|8)$. \end{Statement}

\ssbegin[$\fab(3)$ in the grading $N=8$, where only $\fg_\ev$ is left]{Statement}[$\fab(3)$ in the grading $N=8$, where only $\fg_\ev$ is left]\label{SAb8} For the grading $N=8$ of table \eqref{tabAbGr}
$(d=3$, $\sdim\fg_{-3}=(1|1))$, the only semi-simple maximal graded subalgebra is
\be\label{ansAb8}
\footnotesize\begin{tabular}{llll}
$\sdim\fu$&the components of $\fu$&$\fu$&case, certificate\\ \hline
$24|0$&${\scriptstyle -3{:}}1|0,{\scriptstyle -2{:}}3|0,{\scriptstyle -1{:}}3|0,{\scriptstyle +0{:}}10|0,{\scriptstyle +1{:}}3|0,{\scriptstyle +2{:}}3|0,{\scriptstyle +3{:}}1|0$&$\fo(7)\oplus\fsl(2)$&$\fg_\ev$, max\\
\end{tabular}
\ee
Besides it, the sweep carried out over $\Qee$ returns $48$ maximal graded subalgebras that are not
semi-simple, of $\sdim$ equal to $(16|10)$, $(17|9)$, $(17|11)$, $(17|12)$, $(18|12)$,~$(23|8)$. \end{Statement}

\section{Remarks on the answers}\label{ss3excRem}

The three remarks that follow concern, respectively, the value of the parameter $a$ of the
copies of $\fosp_a(4|2)$ that occur in the answers, the meaning of the marks
\textbf{max(b)}, and the branch of the case $3^o$ which the sweep of \S\ref{ss3excSweep}
does not cover.

\ssbegin[The parameter $a$]{Remark}[The parameter $a$]\label{Ra3exc} Every copy of $\fosp_a(4|2)$ that occurs in the list \eqref{unionAg} has $a$ in the $S_3$-orbit of $3$, namely
$a\in\{3,\nfrac13,-4,-\nfrac43,-\nfrac34,-\nfrac14\}$, and every copy that occurs
in the list \eqref{unionAb} has $a$ in the $S_3$-orbit of $2$, namely
$a\in\{2,\nfrac12,-3,-\nfrac13,-\nfrac32,-\nfrac23\}$. These numbers are not
accidental: the even part of such a~$\fu$ is $\fsl(2)^{\oplus3}$, two of whose
summands lie in $\fg(2)$ (resp.~in $\fo(7)$) as a~subalgebra $\fsl(2)\oplus\fsl(2)$ with
roots of unequal length, and the ratio of the lengths is exactly $3$ for $G_2$
(resp.~$2$ for~$\fo(7)$). For such a~subalgebra $\fu\simeq\fosp_a(4|2)$, the value of $a$ is computed intrinsically: $\fu$ is
simple, so its invariant form is unique up to a~factor, and $(b_1:b_2:b_3)$ with
$b_m=B(H_m,H_m)$ for the canonical $H_m=[E_m,F_m]$ of the three $\fsl(2)$'s
determines $(c_1:c_2:c_3)=(b_1^{-1}:b_2^{-1}:b_3^{-1})$, and hence $a$; it is then
computed once again from the Killing form of the ambient restricted to $\fu$, and
the two answers agree in every case. (Observe that the Killing form of $\fosp_a(4|2)$
itself vanishes identically, so the intrinsic Killing form is of no use here.) %\AO{DL: The procedure was checked on the model \eqref{modelOsp} for several values of~$a$.}

In terms of the invariant \eqref{jinv}: for every copy of
$\fosp_a(4|2)$ that occurs in the answers for $\fag(2)$ (it occurs in four of the five gradings, see the list \eqref{unionAg}), we have $j(3)=\nfrac{2197}{144}=\nfrac{13^3}{144}$; for every copy that occurs in the answers for $\fab(3)$ (in seven of the eight gradings, see the list \eqref{unionAb}), we have $j(2)=\nfrac{343}{36}=\nfrac{7^3}{36}$.
\end{Remark}

\ssbegin[On the marks max(b)]{Remark}[On the marks max(b)]\label{Rb3exc} Each answer marked
\textbf{max(b)} in the tables above is of the same superdimension, with the same
graded components and the same identification as an answer marked \textbf{max} in
the same table, and arises from a~$V$ which is a~$G_0$-translate of the
$\ft$-invariant $V$ of the answer marked \textbf{max}: these are conjugate copies, and the mark (b) is
an artifact of the certificate, not a~weaker result --- except for the answers all
of whose occurrences are marked (b), which are the ones marked $(b)$ in the lists \eqref{unionAg} and \eqref{unionAb}. \end{Remark}

\ssbegin[The irreducible branch of the case $3^o$]{Remark}[The irreducible branch of the
case $3^o$]\label{Rirr3exc} In the sweep of \S\ref{ss3excSweep}, the subalgebra $\fu_0$ of
the case $3^o$ runs over the stabilizers $\fst(V)$ of subspaces $V\subset\fg_{-1}$, i.e.,
over the \textit{reducible} $\fu_0$'s. The other branch --- $\fu_0$ almost simple and
irreducible on $\fg_{-1}$ --- is not parametrized by subspaces; it is, however,
subject to the reduction of Remark \ref{Rthin3exc}, and hence can produce
a~maximal $\fu$ only if the first prolongation $(\fg_{-1},\fu_0)_1$ is non-zero,
which forces $\fu_0$ to be one of the algebras in Cartan's list. For the two
contact gradings, this settles the branch completely. \end{Remark}

The following Statement is what makes the lists for the two contact gradings complete:
the sweep of \S\ref{Ssweep} covers only the reducible $\fu_0$'s, and the Statement shows
that the irreducible ones add nothing there.

\ssbegin[The irreducible branch is empty in both contact gradings]{Statement}[The irreducible branch is empty in both contact gradings]\label{Sirr3exc} For the contact grading of $\fag(2)$, the
principally embedded (see \cite{GL}) $\fsl(2)\subset\fg(2)$ is the only proper almost simple subalgebra of
$\fg(2)$ acting irreducibly on $\fg_{-1}=\Pi(\pi_1)$ \textup{(}by Dynkin's
tables, see \cite{Dy, OV}, the maximal subalgebras of $\fg(2)$ are $\fsl(3)$,
$\fsl(2)\oplus\fsl(2)$, the principal $\fsl(2)$ and the two maximal parabolic
subalgebras, and $\pi_1=\id\oplus\id^*\oplus\Cee$ over $\fsl(3)$, whereas
$\pi_1=(\Cee^2\boxtimes\Cee^2)\oplus(\Cee\boxtimes S^2(\Cee^2))$ over
$\fsl(2)\oplus\fsl(2)$\textup{)}, and, computed over $\Qee$,
\[
(\fg_{-1},\fsl(2))_1=(\fg_{-1},\fc(\fsl(2)))_1=0 .
\]
For the contact grading of $\fab(3)$, the principally embedded $\fsl(2)\subset\fo(7)$ is
\textbf{not} irreducible on $\spin_7$ \textup{(}it lies in $\fg(2)\subset\fo(7)$, and
$\spin_7=\pi_1\oplus\Cee$ over $\fg(2)$, so $\spin_7=S^6(\Cee^2)\oplus\Cee$ as an
$\fsl(2)$-module\textup{)}, and again the first prolongation vanishes. Hence, in
both contact gradings the irreducible branch of the case $3^o$ is
empty. \end{Statement}

\ssec{Related results}\label{ssRelated} For the classification of the maximal subalgebras
of the exceptional simple Lie algebras with Cartan matrix over fields of ``good\rq\rq\ 
characteristic $p>0$, see \cite{PS}.

\section{Open problems}\label{ss3excNot}

\ssbegin[The irreducible branch in the $11$ non-contact gradings]{Problem}[The irreducible branch in the $11$ non-contact gradings]\label{Pirr} Settle the irreducible branch of the case $3^o$ for the
$11$ non-contact gradings of tables \eqref{tabAgGr}, \eqref{tabAbGr}. Remark $\ref{Rirr3exc}$ and
Statement $\ref{Sirr3exc}$ settle it for the two contact gradings, where
$\fg_0=\fc(\text{simple})$ and Dynkin's tables suffice. In the other gradings $\fg_0$ is
not of that form, so its almost simple subalgebras acting irreducibly on $\fg_{-1}$ have
to be enumerated by hand. \end{Problem}

\ssbegin[Maximality of the four answers marked max(b)]{Problem}[Maximality of the four answers marked max(b)]\label{Pmaxb} Prove the maximality of the four answers all of
whose occurrences are marked \textbf{max(b)}, i.e., of $\fsl(2)\oplus\fsl(2|1)$
and $\fosp(1|2)\oplus\fosp(1|2)$  in $\fag(2)$, and of
$\fg(2)\oplus\fosp(1|2)$  and
$\fosp(1|2)\oplus\fosp(1|2)\oplus\fsl(2)$  in $\fab(3)$; for them the torus
$\ft'=\fu\cap\ft$ is too small for Lemma $\ref{Lcert3exc}$ to be conclusive.
\end{Problem}

\ssbegin[Are the two gradings of depth $4$ to be included?]{Problem}[Are the two gradings of depth $4$ to be included?]\label{Pdepth4} Are the two gradings of depth $4$ of Statement
$\ref{StGr3exc}$ to be included? They violate the condition \eqref{cond3exc} --- the
$\fg_0$-action on $\fg_{-k}$ is reducible for some $k$ --- but they are W-gradings in the
weaker sense, and the ambients they produce are of interest in their own right.
\end{Problem}

\ssbegin[The answers in terms of odd reflections]{Question}[The answers in terms of odd reflections]\label{Qoddrefl} Can the answers be stated in terms of odd reflections?
Each of the gradings of tables \eqref{tabOspGr}--\eqref{tabAbGr} is given by an element $h$ of
a~maximal toral subalgebra, i.e., by a~subset of a~system of simple roots; the systems of simple
roots of the three exceptions form a~groupoid under odd reflections, see \cite{Se4}, and the answers
should be equivariant with respect to it. \end{Question}

\ssbegin[Does the Zharov--Zolotarev--Onishchik Claim apply here?]{Question}[Does the Zharov--Zolotarev--Onishchik Claim apply here?]\label{QZZO} Does the Claim of \cite{ZZO}, describing when the
stabilizer of a~subspace $V\subset W$, not necessarily a~sub\textbf{super}space, is
a~Lie superalgebra, apply to $\fg_0=\fc(\fg(2))$ and $\fc(\fo(7))$ acting on $\pi_1$
and on~$\spin_7$? For a~discussion of that Claim, see \cite{LSh}. \end{Question}

\ssbegin[Are the three embeddings already in the literature?]{Question}[Are the three embeddings already in the literature?]\label{Qlit} Are the embeddings $\fosp_a(4|2)\subset\fag(2)$ with
$j(a)=\nfrac{13^3}{144}$, $\fsl(3|1)\subset\fag(2)$, and $\fsl(4|1)\subset\fab(3)$
already in the literature? \cite{ShM} lists the maximal
subalgebras of the matrix Lie superalgebras, but the three exceptions are not matrix in
any preferred way, so they are outside its reach. \end{Question}

\ssbegin[The maximal non-graded subalgebras]{Problem}[The maximal non-graded subalgebras]\label{Pnongraded} Describe the maximal \textbf{non-graded} subalgebras
of the three exceptions, i.e., drop the restriction to graded subalgebras imposed in the formulation \eqref{PROB}. Since the
three exceptions are finite-dimensional, this is a~finite problem, unlike its analog for
the infinite-dimensional ambients of \cite{LSh}. \end{Problem}

\para{Disclosures} 
We thank I.~Shchepochkina for raising the problem and for her help with it; her preprint in Stockholm
University Report no.~32/1988-15  is the source of the method used here. No conflict of interest. 
The data used to support the findings of this study are included within the article. The codes are available from O.~Lozhechnyk upon request until they are posted on GitHub.

\end{document}